\documentclass[10pt]{article}
\usepackage[T1]{fontenc}
\usepackage[latin1]{inputenc} 

\usepackage[french, english]{babel}
\usepackage[fleqn]{amsmath}
\usepackage{amssymb,mathrsfs}
\usepackage{enumitem} 
\usepackage{amsthm}
\usepackage{titlesec}
\usepackage[all]{xy}
{\theoremstyle {definition}  }
{\theoremstyle {definition}  }
{\theoremstyle {plain}  }
{\theoremstyle {plain}  }
{\theoremstyle {plain}  }
{\theoremstyle {plain}  }
{\theoremstyle {plain}  }
{\theoremstyle {plain}  }
\usepackage[a4paper]{geometry}
\def\l@paragraph{\@tocline{3}{0pt}{1pc}{9pc}{}}

\newcommand\M[1]{\mathscr{#1}}

\newcommand\N{\mathbb{N}}

\newcommand\K{\mathbb{K}}
\newcommand\red[1]{\text{Red}\left(#1\right)}

\newcommand\im[1]{\text{im}\left(#1\right)}

\newcommand\imm[1]{\emph{im}\left(#1\right)}

\newcommand\F[1]{\underset{#1}{\longrightarrow}}

\newcommand\vect[1]{\overline{#1}}
\newcommand\EV[1]{\left\langle #1\right\rangle}

\newcommand\V[1]{V^{#1}}
\newcommand\id[1]{\text{Id}_{#1}}
\newcommand\idd[1]{\emph{Id}_{#1}}

\newcommand\tens[1]{\text{T}\left(#1\right)}

\newcommand\lm[1]{\text{lm}\left(#1\right)}

\titleformat{\subsubsection}[runin]
{\normalfont\bfseries}
{\thesubsubsection.}{.5em}{}[.]
\titlespacing{\subparagraph}
{\parindent}{1.5ex plus .1ex minus .2ex}{1.5pt}

\begin{document}

\title{Confluence algebras and acyclicity of the Koszul complex}
\author{CYRILLE CHENAVIER}
\date{}
\maketitle
\paragraph{Abstract.}

The $N$-Koszul algebras are $N$-homogeneous algebras which satisfy an homological property. These algebras are characterised by their Koszul complex: an $N$-homogeneous algebra is $N$-Koszul if and only if its Koszul complex is acyclic. Methods based on computational approaches were used to prove $N$-Koszulness: an algebra admitting a side-confluent presentation is $N$-Koszul if and only if the extra-condition holds. However, in general, these methods do not provide an explicit contracting homotopy for the Koszul complex. In this article we present a way to construct such a contracting homotopy. The property of side-confluence enables us to define specific representations of confluence algebras. These representations provide a candidate for the contracting homotopy. When the extra-condition holds, it turns out that this candidate works. We explicit our construction on several examples.

\tableofcontents
\pagestyle{plain}

\section{Introduction}

\begin{center}

\textbf{An overview on Koszulness and $N$-Koszulness}
\end{center}

\paragraph{Koszul algebras.}

\emph{Koszul algebras} were defined by Priddy in~\cite{MR0265437} as \emph{quadratic algebras} which satisfy a homological property. A quadratic algebra is a graded associative algebra over a field $\K$ which admits a \emph{quadratic presentation} $\EV{X\mid R}$, that is, $X$ is a set of generators and $R$ a set of quadratic relations. If $A$ is a quadratic algebra, the field $\K$ admits a left and right $A$-module structure induced by the $\K$-linear projection $\varepsilon:A\F{}\K$ which maps any generator of $A$ to $0$. A quadratic algebra $A$ is said to be Koszul if the Tor groups $\text{Tor}^A_{n,(m)}\left(\K,\K\right)$ ($n$ is the homological degree and $m$ is graduation induced by the natural graduation over $A$) vanish for $m\neq n$. 

A property of Koszul algebras is that the ground field $\K$ admits a \emph{Koszul resolution}. The name of this resolution is due to the fact that it is inspired by ideas of Koszul (see~\cite{MR0036511}). The \emph{Koszul complex} of a quadratic algebra $A$ which admits a quadratic presentation $\EV{X\mid R}$ is the complex of free left $A$-modules:
\[\cdots\overset{\partial_{n+1}}{\F{}}A\otimes J_n\overset{\partial_n}{\F{}}A\otimes J_{n-1}\F{}\cdots\overset{\partial_4}{\F{}}A\otimes J_3\overset{\partial_3}{\F{}} A\otimes\overline{R} \overset{\partial_2}{\F{}}A\otimes\K X\overset{\partial_1}{\F{}}A\overset{\varepsilon}{\F{}}\K\F{}0,\]
where $\K X$ and $\overline{R}$ denote the vector space spanned by $X$ and the sub vector space of $\K X^{\otimes 2}$ spanned by $R$, respectively, and for every integer $n$ such that $n\geq 2$, we have:
\[J_n=\bigcap_{i=0}^{n-2}\K X^{\otimes i}\otimes\overline{R}\otimes\K X^{\otimes n-N-i}.\]
The differentials of the Koszul complex are defined by the inclusions of $\overline{R}$ in $A\otimes\K X$, of $J_3$ in $A\otimes\overline{R}$ and of $J_n$ in $A\otimes J_{n-1}$ for $n\geq 4$. Then, a quadratic algebra is Koszul if and only if its Koszul complex is acyclic, that is, if and only if the Koszul complex of $A$ is a resolution of $\K$. 

Another characterisation of Koszulness was given by Backelin in~\cite{MR789425} (see also Theorem 4.1 in~\cite[chapter 2]{MR2177131}): a quadratic algebra is Koszul if and only if it is distributive (which means that some lattices defined with $X$ and $R$ are distributive). Moreover, Koszul algebras have been studied through computational approaches based a monomial order, that is, a well founded total order on the set of monomials. In~\cite{MR846601}, Anick used Gr\"obner basis to construct a free resolution of $\K$. This resolution enables us to conclude that an algebra which admits a quadratic Gr\"obner basis is Koszul. In~\cite{MR1608711}, Berger studied quadratic algebras with a \emph{side-confluent presentation}~\footnote{This notion corresponds to the one of \emph{X-confluent algebra} in~\cite{MR1608711}. However, we prefer to use our terminology because the property of confluence depends on the presentation.}. The latter is a transcription of the notion of quadratic Gr\"obner basis using some linear operators. More precisely, we can associate with any quadratic presentation $\EV{X\mid R}$ of $A$ an unique linear projector $S$ of $\K X\otimes\K X$. This projector maps any element of $\K X\otimes\K X$ to a better one with respect to the monomial order. The presentation $\EV{X\mid R}$ is said to be side-confluent if there exists an integer $k$ such that:
\[\EV{S\otimes\id{\K X},\id{\K X}\otimes S}^k=\EV{\id{\K X}\otimes S,S\otimes\id{\K X}}^k,\]
where $\EV{t,s}$ denotes the product $\cdots sts$ with $k$ factors. The algebra presented by:
\[\EV{s_1,s_2\mid\EV{s_1,s_2}^k=\EV{s_2,s_1}^k,\ s_i^2=s_i,\ i=1,2},\]
is naturally associated with a side-confluent presentation. This algebra is the \emph{confluence algebra of degree k.} In~\cite[Section 5]{MR1608711}, Berger used specific representations of these algebras to construct a contracting homotopy for the Koszul complex of an algebra which admits a side-confluent presentation. This construction enables us to conclude that a quadratic algebra admitting a side-confluent presentation is Koszul.

\paragraph{$N$-Koszul algebras.}

Let $N$ be an integer such that $N\geq 2$. An \emph{N-homogeneous algebra} is a graded associative algebra over a field $\K$ which admits an \emph{N-homogeneous presentation} $\EV{X\mid R}$, that is, $X$ is a set of generators and $R$ is a set of $N$-homogeneous relations. In~\cite{MR1832913} the notion of Koszul algebra is extended to the notion of \emph{N-Koszul algebra.} An $N$-homogeneous algebra $A$ is said to be \emph{N}-Koszul if the Tor groups $\text{Tor}^A_{n,(m)}\left(\K,\K\right)$ vanish for $m\neq l_N(n)$, where $l_N$ is the function defined by:
\[l_N(n)=\left\{
\begin{split}
&kN,\ \ \ \ \ \ \text{if}\ n=2k,\\
&kN+1,\ \text{if}\ n=2k+1.
\end{split}
\right.\]
We remark that a $2$-Koszul algebra is precisely a Koszul algebra. Thus, the notion of $N$-Koszul algebra generalises the one of Koszul algebra.

In the same paper, Berger defined the \emph{Koszul complex} of an $N$-homogeneous algebra. Let $\EV{X\mid R}$ be an $N$-homogeneous presentation of $A$. The Koszul complex of $A$ is the complex of left $A$-modules:
\[\cdots\overset{\partial_{n+1}}{\F{}}A\otimes J^N_n\overset{\partial_n}{\F{}}A\otimes J^N_{n-1}\F{}\cdots\overset{\partial_4}{\F{}}A\otimes J^N_3\overset{\partial_3}{\F{}} A\otimes\overline{R} \overset{\partial_2}{\F{}}A\otimes\K X\overset{\partial_1}{\F{}}A\overset{\varepsilon}{\F{}}\K\F{}0,\]
where the vector spaces $J^N_n$ are defined by:
\[J^N_n=\bigcap_{i=0}^{l_N(n)-N}\K X^{\otimes i}\otimes\overline{R}\otimes\K X^{\otimes l_N(n)-N-i}.\]

As in the quadratic case, this complex characterises the property of $N$-Koszulness: an $N$-homogeneous algebra is $N$-Koszul if and only if its Koszul complex is acyclic (see~\cite[Proposition 2.12]{MR1832913}). This complex also find applications in the study of higher Koszul duality (see~\cite{MR3110804}).

Berger studied the property of $N$-Koszulness using monomial orders. As in the quadratic case, there exists a unique linear projector $S$ of $\K X^{\otimes N}$ associated with an $N$-homogeneous presentation of $A$ which maps any element to a better one with respect to the monomial order. Then, a presentation is \emph{side-confluent} if for every integer $m$ such that $N+1\leq m\leq 2N-1$, there exists an integer $k$ which satisfies:
\[\EV{S\otimes\id{\K X^{\otimes m-N}},\id{\K X^{\otimes m-N}}\otimes S}^k=\EV{\id{\K X^{\otimes m-N}}\otimes S,S\otimes\id{\K X^{\otimes m-N}}}^k.\]
Contrary to the quadratic case, an algebra admitting a side-confluent presentation is not necessarily $N$-Koszul. Indeed, when the set $X$ is finite, such an algebra is $N$-Koszul if and only if the \emph{extra-condition} holds (see \cite[Proposition 3.4]{MR1832913}). The extra-condition is stated as follows:
\[(ec): \ \ \ \ \left(\K X^{\otimes m}\otimes\overline{R}\right)\cap\left(\overline{R}\otimes\K X^{\otimes m}\right)\ \subset\ \K X^{\otimes m-1}\otimes\overline{R}\otimes\K X,\ \text{for every}\ 2\leq m\leq N-1,\]
We group these hypothesis in the following definition:

\begin{quote}

\textbf{\ref{Definition of extra-confluent presentation} Definition.} Let $A$ be an $N$-homogeneous algebra. A side-presentation $\EV{X\mid R}$ such that $X$ is finite and the extra-condition holds is said to be \emph{extra-confluent.}

\end{quote}

\paragraph{Our problematic.} We deduce of the works from~\cite{MR1832913} that the Koszul complex of an algebra $A$ admitting an extra-confluent presentation is acyclic. However, there does not exist an explicit contracting homotopy for the Koszul complex of $A$. The purpose of this paper is to construct such a contracting homotopy. For the quadratic case, our contracting homotopy is the one constructed in~\cite[Section 5]{MR1608711}.

\begin{center}

\textbf{Our results}

\end{center}
We present the different steps of our construction. Recall that an extra-confluent presentation needs a monomial order. Thus, in what follows, we work with a monomial order. For every integer $m$, we denote by $X^{(m)}$ the set of words of length $m$.

\paragraph{Reduction pairs associated with a presentation.} In the way to construct our contracting homotopy, we will associate with any $N$-homogeneous presentation $\EV{X\mid R}$ such that $X$ is finite a family $P_{n,m}=\left(F_1^{n,m},F_2^{n,m}\right)$, where $F_1^{n,m}$ and $F_2^{n,m}$ are linear projectors of $\K X^{(m)}$. The pair $P_{n,m}$ is called the \emph{reduction pair of bi-degree} $(n,m)$ associated with $\EV{X\mid R}$. We point the fact that the finiteness condition over $X$ will be necessary to define the operators $F_i^{n,m}$. Moreover, these operators satisfy the following condition: for any $w\in X^{(m)}$, $F_i^{n,m}(w)$ is either equal to $w$ or is a sum of monomials which are strictly smaller than $w$ with respect to the monomial order. The linear projectors of $\K X^{(m)}$ satisfying the previous condition are called \emph{reduction operators relatively to} $X^{(m)}$. The set of reduction operators relatively to $X^{(m)}$ admits a lattice structure (we will recall it in Section~\ref{Reduction operators confluence algebras}). This structure plays an essential role in our constructions. A pair $\left(T_1,T_2\right)$ of reduction operators relatively to $X^{(m)}$ is said to be \emph{confluent} if there exists an integer $k$ such that we have the following equality in $\text{End}\left(\K X^{(m)}\right)$:
\[\EV{T_1,T_2}^k=\EV{T_2,T_1}^k.\]
Then, our first result is:

\begin{quote}

\textbf{\ref{Confluence} Theorem} \textit{Let $A$ be an $N$-homogeneous algebra admitting a side-confluent presentation $\EV{X\mid R}$, where $X$ is a finite set. The reduction pairs associated with $\EV{X\mid R}$ are confluent.}

\end{quote}

\paragraph{The left bound.}

The reduction pairs associated with a side-confluent presentation $\EV{X\mid R}$ enable us to define a family of representations of confluence algebras in the following way:
\[\begin{split}
  \varphi^{P_{n,m}}\colon\EV{s_1,s_2\mid\EV{s_1,s_2}^{k_{n,m}}=\EV{s_2,s_1}^{k_{n,m}},\ s_i^2=s_i,\ i=1,2}&\longrightarrow\text{End}\left(\K X^{(m)}\right),\\
  s_i&\longmapsto F_i^{n,m}
\end{split}
\]
where the integer $k_{n,m}$ satisfies:
\[\EV{F_1^{n,m},F_2^{n,m}}^{k_{n,m}}=\EV{F_2^{n,m},F_1^{n,m}}^{k_{n,m}}.\]
For every integers $n$ and $m$ we will consider a specific element in $\M{A}_{k_{n,m}}$:
\[\gamma_1=(1-s_2)\left(s_1+s_1s_2s_1+\cdots +\EV{s_2,s_1}^{2i+1}\right),\]
where the integer $i$ depends on $k_{n,m}$. The shape of this element will be motivated in Section~\ref{The contracting homotopy in small degree}. In Section~\ref{Construction} we will use the elements $\varphi^{P_{n,m}}\left(\gamma_1\right)$ to construct a family of $\K$-linear maps
\[\begin{split}
&h_0:A\F{} A\otimes\K X,\\
&h_1:A\otimes\K X\F{} A\otimes\overline{R},\\
&h_2:A\otimes\overline{R}\F{}A\otimes J^N_3,\\
&h_n:A\otimes J^N_n\F{} A\otimes J^N_{n+1},\ \text{for}\ n\geq 3.
\end{split}\]
The family $\left(h_n\right)_n$ is called the \emph{left bound of} $\EV{X\mid R}$. In Proposition~\ref{relation r}, we will show that the left bound of $\EV{X\mid R}$ is a contracting homotopy for the Koszul complex of the algebra presented by $\EV{X\mid R}$ if and only if $\EV{X\mid R}$ satisfies some identities, called the \emph{reduction relations}.

\paragraph{Extra-confluent presentations and reduction relations.}

Finally we will show that $(ec)$ implies that the reduction relations hold. Then, our main result is stated as follows:

\begin{quote}

\textbf{\ref{Theorem} Theorem} \textit{Let $A$ be an $N$-homogeneous algebra. If $A$ admits an extra-confluent presentation $\EV{X\mid R}$, then the left bound of $\EV{X\mid R}$ is a contracting homotopy for the Koszul complex of $A$.}

\end{quote}

\begin{center}

\textbf{Organisation}

\end{center}
In Section~\ref{Preliminaries} we recall how we can construct the Koszul complex of an $N$-homogeneous algebra. We also recall the definition of an extra-confluent presentation. In Section~\ref{The contracting homotopy in small degree} we make explicit our construction in small homological degree. In Section~\ref{Reduction operators confluence algebras} we recall the definitions of confluence algebras and reduction operators. We also recall the link between reduction operators and representations of confluence algebras. In Section~\ref{The left bound of a side-confluent presentation} we construct the contracting homotopy in terms of confluence. As an illustration of our construction we provide in Section~\ref{Examples} three examples: the \emph{symmetric algebra}, \emph{monomial algebras} which satisfy the \emph{overlap properties} and the \emph{Yang-Mills algebra over two generators.}\\\\
\textbf{Acknowledgement.} The author wish to thank Roland Berger for helpful discussions. This work is supported by the Sorbonne-Paris-Cit\'e IDEX grant Focal and the ANR grant ANR-13-BS02-0005-02 CATHRE. 

\section{Preliminaries}\label{Preliminaries} 

\subsection{The Koszul complex}\label{The Koszul complex}

\subsubsection{Conventions and notation}\label{Conventions and notations}

We denote by $\K$ a field. We say vector space and algebra instead of $\K$-vector space and $\K$-algebra, respectively. We consider only associative algebras. Given a set $X$, we denote by $\EV{X}$ and $\K X$ the free monoid and the vector space spanned by $X$, respectively. For every integer $m$, we denote by $X^{(m)}$ the subset of $\EV{X}$ of words of length $m$.

We write $V=\K X$. We identify $\K X^{(m)}$ and the free algebra $\K\EV{X}$ spanned by $X$ to $\V{\otimes m}$ and to the tensor algebra $\tens{V}$ over $V$, respectively.

Let $A$ be an algebra. A \emph{presentation} of $A$ is a pair $\EV{X\mid R}$, where $X$ is a set and $R$ is a subset $\K\EV{X}$ such that $A$ is isomorphic to the quotient of $\K\EV{X}$ by $I(R)$, where $I(R)$ is the two-sided ideal of $\K\EV{X}$ spanned by $R$. The isomorphism from $A$ to $\K\EV{X}/I(R)$ is denoted by $\psi_{\EV{X\mid R}}$. For every $f\in\K\EV{X}$, we denote by $\overline{f}$ the image of $f$ through the canonical projection of $\K\EV{X}$ over $A$.

Let $N$ be an integer such that $N\geq 2$. An $N$-\emph{homogeneous presentation} of $A$ is a presentation $\EV{X\mid R}$ of $A$ such that $R$ is included in $\K X^{(N)}$. In this case, the two-sided ideal $I(R)$ is the direct sum of vector spaces $I(R)_m$ defined by $I(R)_m=0$ if $m<N$, and
\[I(R)_m=\sum_{i=0}^{m-N}\V{\otimes i}\otimes\overline{R}\otimes\V{\otimes m-N-i}\ \text{if}\ m\geq N,\]
where $\overline{R}$ denotes the subspace of $\V{\otimes N}$ spanned by $R$. An $N$-\emph{homogeneous algebra} is a graded algebra $A=\bigoplus_{m\in\N}A_m$ which admits an $N$-homogeneous presentation $\EV{X\mid R}$ such that for every integer $m$, $\psi_{\EV{X\mid R}}$ induces a $\K$-linear isomorphism from $A_m$ to $\V{\otimes m}/I(R)_m$:
\[\begin{split}
A&=\bigoplus_{m\in\N}A_m\\
&\simeq\K\oplus\V{}\oplus\cdots\oplus\V{\otimes N-1}\oplus\frac{\V{\otimes N}}{\overline{R}}\oplus\frac{\V{\otimes N+1}}{V\otimes\overline{R} + \overline{R}\otimes V}\oplus\cdots
\end{split}\]
We denote by $\varepsilon:A\F{}\K$ the projection which maps $1_A$ to $1_\K$ and $A_m$ to $0$ for every $m\geq 1$.

\subsubsection{The construction of the Koszul complex}\label{The construction of the Koszul complex}

Let $A$ be an $N$-homogeneous algebra and let $\EV{X\mid R}$ be an $N$-homogeneous presentation of $A$. We write $V=\K X$. We consider the family of vector spaces $(J^N_n)_n$ defined by $J^N_0=\K$, $J^N_1=V$, $J^N_2=\overline{R}$ and for every integer $n\geq 3$
\[J^N_n=\bigcap_{i=0}^{l_N(n)-N}\V{\otimes i}\otimes\overline{R}\otimes\V{\otimes l_N(n)-N-i},\]
where the function $l_N:\N\F{}\N$ is defined by
\[l_N(n)=\left\{
\begin{split}
&kN,\ \ \ \ \ \ \text{if}\ n=2k,\\
&kN+1,\ \text{if}\ n=2k+1.
\end{split}
\right.\]
When there is no ambiguity, we write $J_n$ instead of $J^N_n$.

Let $n$ be an integer. For every $w\in X^{\left(l_N(n+1)\right)}$, let $w_1\in X^{\left(l_N(n+1)-l_N(n)\right)}$ and $w_2\in X^{\left(l_N(n)\right)}$ such that $w=w_1w_2$. Let us consider the $A$-linear map
\[
\begin{split}
F_{n+1}\colon A\otimes\V{\otimes l_N(n+1)}&\F{} A\otimes\V{\otimes l_N(n)}.\\
1_A\otimes w&\longmapsto \overline{w_1}\otimes w_2
\end{split}
\]
Recall from \cite[Section 3]{MR1832913} that the \emph{Koszul complex} of $A$ is the complex $\left(K_\bullet,\partial\right)$
\[\cdots\overset{\partial_{n+1}}{\F{}}A\otimes J_n\overset{\partial_n}{\F{}}A\otimes J_{n-1}\F{}\cdots\overset{\partial_2}{\F{}}A\otimes J_1\overset{\partial_1}{\F{}}A\overset{\varepsilon}{\F{}}\K\F{}0,\]
where $\partial_{n}$ is the restriction of $F_{n}$ to $A\otimes J_n$. In particular, the map $\partial_1$ is defined by $\partial_1\left(1_A\otimes x\right)~=~\overline{x}$ for every $x\in X$.

\subsubsection{Remark}
The two following remarks show that the Koszul complex is well-defined:

\begin{enumerate}
\item Let $n$ be an integer. The vector space $J_{n+1}$ is included in $\V{\otimes l_N(n+1)-l_N(n)}\otimes J_n$. Thus, the vector space $F_{n+1}\left(A\otimes J_{n+1}\right)$ is included in $A\otimes J_n$.
\item Let $n$ be an integer such that $n\geq 1$. The vector space $J_{n+1}$ is included in $R\otimes J_{n-1}$. Thus, the restriction of $F_nF_{n+1}$ to $A\otimes J_{n+1}$ vanishes.
\end{enumerate}

\subsubsection{Example}\label{Koszul complex of the Yang-Mills algebra}

We consider the example from~\cite[Section 6.3]{MR3299599}. Let $A$ be the \emph{Yang-Mills algebra over 2 generators}: this algebra is the $3$-homogeneous algebra presented by 
\[\EV{x_1,x_2\mid x_2x_1x_1-2x_1x_2x_1+x_1x_1x_2,\ x_2x_2x_1-2x_2x_1x_2+x_1x_2x_2}.\]
The map $\partial_2:A\otimes\overline{R}\F{} A\otimes V$ is defined by
\[\begin{split}
\partial_2\left(1_A\otimes x_2x_1x_1-2x_1x_2x_1+x_1x_1x_2\right)&=\overline{x_2x_1}\otimes x_1-2\overline{x_1x_2}\otimes x_1+\overline{x_1x_1}\otimes x_2,\ \text{and}\\
\partial_2\left(1_A\otimes x_2x_2x_1-2x_2x_1x_2+x_1x_2x_2\right)&=\overline{x_2x_2}\otimes x_1-2\overline{x_2x_1}\otimes x_2+\overline{x_1x_2}\otimes x_2.
\end{split}\]
The vector space $J_3=\left(V\otimes\overline{R}\right)\cap\left(\overline{R}\otimes V\right)$ is the one-dimensional vector space spanned by
\[\begin{split}
v&=x_2(x_2x_1x_1-2x_1x_2x_1+x_1x_1x_2)+x_1(x_2x_2x_1-2x_2x_1x_2+x_1x_2x_2)\\
&=(x_2x_2x_1-2x_2x_1x_2+x_1x_2x_2)x_1+(x_2x_1x_1-2x_1x_2x_1+x_1x_1x_2)x_2.
\end{split}\]
The map $\partial_3:A\otimes J_3\F{} A\otimes\overline{R}$ is defined by
\[\partial_3\left(1_A\otimes v\right)=\overline{x_2}\otimes\left(x_2x_1x_1-2x_1x_2x_1+x_1x_1x_2\right)+\overline{x_1}\otimes\left(x_2x_2x_1-2x_2x_1x_2+x_1x_2x_2\right).\]

\subsection{Side-confluent presentations}\label{Side-confluent presentations}

Through this section we fix an $N$-homogeneous algebra $A$ and an $N$-homogeneous presentation $\EV{X\mid R}$ of $A$. We assume that $X$ is a totally ordered set. We write $V=\K X$.

\subsubsection{Reductions}\label{Reductions}

For every integer $m$, the set $X^{(m)}$ is totally ordered for the lexicographic order induced by the order over $X$. For every $f\in\V{\otimes m}\setminus\{0\}$, we denote by $\lm{f}$ the greatest element of $X^{(m)}$ occurring in the decomposition of $f$. We denote by $\text{lc}(f)$ the coefficient of $\lm{f}$ in the decomposition of $f$. Let
\[R'=\left\{\frac{1}{\text{lc}\left(f\right)}f,\ f\in R\right\}.\]
Then, $\EV{X\mid R'}$ is an $N$-homogeneous presentation of $A$. Thus, we can assume that $\text{lc}\left(f\right)$ is equal to~1 for every $f\in R$. 

For every $w_1,w_2\in\EV{X}$ and every $f\in R$, let $r_{w_1fw_2}$ be the $\K$-linear endomorphism of $\tens{V}$ defined on the basis $\EV{X}$ in the following way:
\[r_{w_1fw_2}\left(w\right)=\left\{
\begin{split}
&w_1\left(\lm{f}-f\right)w_2,\ \text{if}\ w=w_1\lm{f}w_2,\\
&w,\ \text{otherwise}.
\end{split}
\right.\]
Taking the terminology of~\cite{MR506890}, the morphisms $r_{w_1fw_2}$ are called the \emph{reductions of} $\EV{X\mid R}$.

\subsubsection{Normal forms}\label{Normal forms}

An element $f\in\tens{V}$ is said to be a \emph{normal form for} $\EV{X\mid R}$ if $r(f)=f$ for every reduction $r$ of $\EV{X\mid R}$. Given an element $f$ of $\tens{V}$, a \emph{normal form of} $f$ is a normal form $g$ such that there exist reductions $r_1,\cdots ,r_n$ satisfying $g=r_1\cdots r_n(f)$. In this case, we have $\overline{f}=\overline{g}$.

The presentation $\EV{X\mid R}$ is said to be \emph{reduced} if, for every $f\in R$, $\lm{f}-f$ is a normal form for $\EV{X\mid R}$ and $\lm{f}$ is a normal form for $\EV{X\mid R\setminus\{f\}}$. From this moment, all the presentations are assumed to be reduced.

\subsubsection{Critical branching}\label{Critical branching}

A \emph{critical branching of} $\EV{X\mid R}$ is a 5-tuple $(w_1,w_2,w_3,f,g)$ where $f,g~\in ~ R$ and $w_1,w_2,w_3$ are non empty words such that:
\[\begin{split}
w_1w_2&=\lm{f},\ \text{and}\\
w_2w_3&=\lm{g}.
\end{split}\]
The word $w_1w_2w_3$ is the \emph{source} of this critical branching.

\subsubsection{The operator of a presentation}\label{The operator of a presentation}

Let $S$ be the endomorphism of $\V{\otimes N}$ defined on the basis $X^{(N)}$ in the following way:
\[S\left(w\right)=\left\{
\begin{split}
&\lm{f}-f,\ \text{if there exists}\ f\in R\ \text{such that}\ w=\lm{f}, \\
&w,\ \text{otherwise}.
\end{split}
\right.\]
The operator $S$ is \emph{the operator of} $\EV{X\mid R}$. The presentation $\EV{X\mid R}$ is reduced. Thus, $S$ is well-defined and is a projector. The kernel of $S$ is equal to $\overline{R}$. If $w\in X^{(N)}$ is a normal form, then $S(w)$ is equal to $w$. If $w$ is not a normal form, then $S(w)$ is strictly smaller than $w$.

\subsubsection{Definition}\label{Definition of a side-confluent presentation}

The presentation $\EV{X\mid R}$ is said to be \emph{side-confluent} if for every integer $m$ such that $1\leq m\leq N-1$, there exists an integer $k$ such that:
\[\EV{\id{\V{\otimes m}}\otimes S,S\otimes\id{\V{\otimes m}}}^k=\EV{S\otimes\id{\V{\otimes m}},\id{\V{\otimes m}}\otimes S}^k,\]
where $\EV{t,s}^k$ denotes the product $\cdots sts$ with $k$ factors.\\

The Diamond Lemma (~\cite[Theorem 1.2]{MR506890}) implies the following:

\subsubsection{Proposition}\label{Diamond lemma}

\textit{Let $A$ be an $N$-homogeneous algebra. Assume that $A$ admits a side-confluent presentation $\EV{X\mid R}$.} Then, the following hold:
\begin{enumerate}
\item \textit{Every element of} $\tens{V}$ \textit{admits a unique normal form for} $\EV{X\mid R}$.
\item \textit{The set $\left\{\overline{w},\ w\in\EV{X}\ \text{is a normal form}\right\}$ is a basis of A.} 
\item \textit{An element of} $\tens{V}$ belongs to $I(R)$ \textit{if and only if its normal form is equal to} 0.
\end{enumerate}

\begin{proof}

Let $S$ be the operator of $\EV{X\mid R}$. Let $(w_1,w_2,w_3,f,g)$ be a critical branching of $\EV{X\mid R}$. Let $m$ be the length of $w=w_1w_2w_3$. The presentation $\EV{X\mid R}$ being $N$-homogeneous, we have $N+1\leq m\leq 2N-1$. Thus, there exists an integer $k$ such that:
\[\EV{\id{\V{\otimes m -N}}\otimes S,S\otimes\id{\V{\otimes m -N}}}^k(w)=\EV{S\otimes\id{\V{\otimes m -N}},\id{\V{\otimes m -N}}\otimes S}^k(w).\]
Thus, there exist two sequences of reductions $r_1,\cdots ,r_n$ and $r'_1,\cdots ,r'_l$ such that $r_1\cdots r_n\left(\left(\lm{f}-f\right)w_3\right)$ is equal to $r'_1\cdots r'_l\left(w_1\left(\lm{g}-g\right)\right)$.  We deduce from~\cite[Theorem 1.2]{MR506890} that every element $f\in\tens{V}$ admits a unique normal form for $\EV{X\mid R}$ and that $\left\{\overline{w},\ w\in\EV{X}\ \text{is a normal form}\right\}$ is a basis of $A$. Thus, the two first points hold.

Let us show the third point. Let $f$ be an element of $\tens{V}$ and let $\widehat{f}$ be its unique normal form. We write:
\[\widehat{f}=\sum_{i\in I}\lambda_iw_i,\]
where $w_i\in\EV{X}$ are normal forms. Then, $\overline{f}$ is equal to $\sum_{i\in I}\lambda_i\overline{w_i}$. From the second point, $\overline{f}$ is equal to $0$ if and only if $\lambda_i$ is equal to $0$ for every $i\in I$.

\end{proof}

\subsubsection{Lemma}\label{Lambda}

\textit{Assume that the presentation} $\EV{X\mid R}$ \textit{is side-confluent. Let S be the operator of} $\EV{X\mid R}$. \textit{For every integer m such that} $N+1\leq m\leq 2N-1$, \textit{there exists an integer k such that:}
\[\begin{split}
&\EV{\id{\V{\otimes m}}-\id{\V{\otimes m -N}}\otimes S,\id{\V{\otimes m}}-S\otimes\id{\V{\otimes m -N}}}^k\\
=&\EV{\id{\V{\otimes m}}-S\otimes\id{\V{\otimes m -N}},\id{\V{\otimes m}}-\id{\V{\otimes m -N}}\otimes S}^k.
\end{split}\]
\textit{Moreover, for every} $w\in X^{(m)}$ \textit{such that} $\id{\V{\otimes m -N}}\otimes S(w)$ \textit{and} $S\otimes\id{\V{\otimes m -N}}(w)$ \textit{are different from w, we have}:
\[\lm{\left(\EV{\id{\V{\otimes m}}-\id{\V{\otimes m -N}}}\otimes S,\id{\V{\otimes m}}-S\otimes\id{\V{\otimes m -N}}}^k(w)\right)=w.\]

\begin{proof}
We write $S_1~=~\id{\V{\otimes m -N}}\otimes S$ and $S_2~=~S\otimes\id{\V{\otimes m -N}}$.\\

The presentation $\EV{X\mid R}$ is side-confluent. Thus, there exists $k\in\N$ such that $\EV{S_2,S_1}^k$ is equal to $\EV{S_1,S_2}^k$. The morphisms $S_1$ and $S_2$ being projectors we show by induction that for every integer $j$ we have:
\[\begin{split}
\EV{\id{\V{\otimes m}}-S_1,\id{\V{\otimes m}}-S_2}^j&=\id{\V{\otimes m}}+\sum_{i=1}^{j-1}(-1)^i\left(\EV{S_1,S_2}^i+\EV{S_2,S_1}^i\right)+(-1)^{j}\EV{S_1,S_2}^j,\\
\EV{\id{\V{\otimes m}}-S_2,\id{\V{\otimes m}}-S_1}^j&=\id{\V{\otimes m}}+\sum_{i=1}^{j-1}(-1)^i\left(\EV{S_1,S_2}^i+\EV{S_2,S_1}^i\right)+(-1)^{j}\EV{S_2,S_1}^j.
\end{split}\]
In particular we have:
\[\EV{\id{\V{\otimes m}}-S_2,\id{\V{\otimes m}}-S_1}^k=\EV{\id{\V{\otimes m}}-S_1,\id{\V{\otimes m}}-S_2}^k.\]

Moreover, if $w\in X^{(m)}$ is such that $S_1(w)$ and $S_2(w)$ are different from $w$, then $S_1(w)$ and $S_2(w)$ are strictly smaller than $w$. We deduce from the relation 
\[\EV{\id{\V{\otimes m}}-S_1,\id{\V{\otimes m}}-S_2}^k(w)=w+\sum_{i=1}^{k-1}(-1)^i\left(\EV{S_1,S_2}^i+\EV{S_2,S_1}^i\right)(w)+(-1)^{k}\EV{S_1,S_2}^k(w),\]
that $\lm{\EV{\id{\V{\otimes m}}-S_1,\id{\V{\otimes m}}-S_2}^k(w)}$ is equal to $w$.
\end{proof}

\subsubsection{Example}\label{Yang-Mills side-confluent}

We consider the presentation from Example~\ref{Koszul complex of the Yang-Mills algebra} of the Yang-Mills algebra over two generators with the order $x_1<x_2$. The operator $S\in\text{End}\left(\V{\otimes 3}\right)$ of this presentation is defined on the basis $X^{(3)}$ by
\[S(w)=\left\{\begin{split}
&2x_1x_2x_1-x_1x_1x_2,\ \text{if}\ w=x_2x_1x_1, \\
&2x_2x_1x_2-x_1x_2x_2,\ \text{if}\ w=x_2x_2x_1, \\
&w,\ \text{otherwise}.
\end{split}\right.\]
This presentation admits exactly one critical branching:
\[\left(x_2,x_2x_1,x_1,x_2x_1x_1-2x_1x_2x_1+x_1x_1x_2,\ x_2x_2x_1-2x_2x_1x_2+x_1x_2x_2\right).\]
We have:
\[\begin{split}
\EV{S\otimes\id{V},\id{V}\otimes S}^2\left(x_2x_2x_1x_1\right)&=\EV{\id{V}\otimes S,S\otimes\id{V}}^2\left(x_2x_2x_1x_1\right)\\
&=x_2x_1x_2x_1 - 2x_1x_2x_1x_2 + x_1x_1x_2x_2.
\end{split}\]
Moreover, for every $w\in X^{(4)}$ which is different from $x_2x_2x_1x_1$, we check that $\EV{S\otimes\id{V},\id{V}\otimes S}^2\left(w\right)$ is equal to $\EV{\id{V}\otimes S,S\otimes\id{V}}^2\left(w\right)$. Thus, we have:
\[\EV{S\otimes\id{V},\id{V}\otimes S}^2=\EV{\id{V}\otimes S,S\otimes\id{V}}^2.\]

For every $w\in X^{(5)}$ we check that $\EV{S\otimes\id{\V{\otimes 2}},\id{\V{\otimes 2}}\otimes S}^2(w)$ and $\EV{\id{\V{\otimes 2}}\otimes S,S\otimes\id{\V{\otimes 2}}}^2(w)$ are equal. Thus, we have: \[\EV{S\otimes\id{\V{\otimes 2}},\id{\V{\otimes 2}}\otimes S}^2=\EV{\id{\V{\otimes 2}}\otimes S,S\otimes\id{\V{\otimes 2}}}^2.\]
We conclude that the presentation from Example~\ref{Koszul complex of the Yang-Mills algebra} with the order $x_1<x_2$ is side-confluent.

\subsection{Extra-confluent presentations}\label{Extra-confluent presentations}

\subsubsection{The extra-condition}\label{The extra-condition}

Let $A$ be an $N$-homogeneous algebra. Assume that $A$ admits a side-confluent presentation $\EV{X\mid R}$ where $X$ is a totally ordered finite set. Recall from~\cite[Section 3]{MR1832913} that the Koszul complex of $A$ is acyclic if and only if the \emph{extra-condition} holds. The extra-condition is stated as follows:
\[\left(\K X^{(n)}\otimes\overline{R}\right)\cap\left(\overline{R}\otimes\K X^{(n)}\right)\ \subset\ \K X^{(n-1)}\otimes\overline{R}\otimes\K X,\ \text{for every}\ 2\leq n\leq N-1.\]

\subsubsection{Definition}\label{Definition of extra-confluent presentation}

Let $A$ be an $N$-homogeneous algebra. A side-presentation $\EV{X\mid R}$ such that $X$ is finite and the extra-condition holds is said to be \emph{extra-confluent.}

\subsubsection{Remark}\label{Extra-confluent presentation for quadratic algebras}

If $N=2$, the extra-condition is an empty condition. Thus, in this case, the notions of extra-confluent presentation and side-confluent presentation coincide.\\

An extra-confluent presentation has the following interpretation in terms of critical branching:

\subsubsection{Proposition}\label{extra condition and critical pair}

\textit{Let $A$ be an $N$-homogeneous algebra. Assume that $A$ admits an extra-confluent presentation $\EV{X\mid R}$. Let $w=x_1\cdots x_m$ be the source of a critical branching of $\EV{X\mid R}$. The word $x_{m-N}\cdots x_{m-1}$ is not a normal form for $\EV{X\mid R}$.}

\begin{proof}
The presentation $\EV{X\mid R}$ is $N$-homogeneous. In particular, we have $N+1\leq m\leq 2N-1$. If $m=N+1$, there is nothing to prove. Thus, we assume that $m$ is greater than $N+2$.

Let $S$ be the operator of $\EV{X\mid R}$. We write
\[S_1=S\otimes\id{\otimes m-N}\ \text{and}\ S_2=\id{\otimes m-N}\otimes S.\]
The presentation $\EV{X\mid R}$ is side-confluent. Thus, from Lemma~\ref{Lambda}, there exists an integer $k$ such that
\[\EV{\id{\V{\otimes m}}-S_2,\id{\V{\otimes m}}-S_1}^k=\EV{\id{\V{\otimes m}}-S_1,\id{\V{\otimes m}}-S_2}^k.\]
We denote by $\Lambda$ this common morphism. By hypothesis, $S_1(w)$ and $S_2(w)$ are different from $w$. From Lemma~\ref{Lambda}, $\lm{\Lambda(w)}$ is equal to $w$.

The image of $\Lambda$ is included in $\im{\id{\K X^{(m)}}-S_1}\cap\im{\id{\K X^{(m)}}-S_2}$ that is, $\ker\left(S_1\right)\cap\ker\left(S_2\right)$. The latter is equal to $\overline{R}\otimes\K X^{(m-N)}\cap\K X^{(m-N)}\otimes\overline{R}$. The presentation $\EV{X\mid R}$ satisfies the extra-condition. Thus, the image of $\Lambda$ is included in $\K X^{(m-N-1)}\otimes\overline{R}\otimes\K X$. In particular, there exist $w_1,\cdots ,w_l\in X^{(m-N-1)}$, $f_1,\cdots , f_l\in R$, $x_1,\cdots ,x_l\in X$ and $\lambda_1,\cdots\lambda_l\in\K$ which satisfy
\[\Lambda(w)=\sum_{i=1}^l\lambda_iw_if_ix_i.\]
Thus, $\lm{\Lambda(w)}=w$ is equal to $w_i\lm{f_i}x_i$ for some $1\leq i\leq l$. We conclude that $x_{n-N}\cdots x_{m-1}$ is equal to $\lm{f_i}$. In particular, it is not a normal form.
\end{proof}

\subsubsection{Remark}

Let $A$ be the algebra presented by $\EV{x<y\mid xyx}$. This presentation is side-confluent. There is only one critical branching: $(xy,x,yx,xyx,xyx)$. The source $xyxyx$ of this critical has length 5. We deduce from Proposition~\ref{extra condition and critical pair} that the extra-condition does not hold.

Let us check that the Koszul complex of $A$ is not acyclic: the vector space $J_3$ is reduced to $\{0\}$ and  the map $\partial_2:A\otimes\overline{R}\F{} A\otimes V$ is defined by $\partial_2\left(1_A\otimes xyx\right)=\overline{xy}\otimes x$. In particular, $\overline{xy}\otimes xyx$ belongs to the kernel of $\partial_2$. Thus, we have a strict inclusion $\im{\partial_3}\ \subsetneqq\ \ker\left(\partial_2\right)$. 

\subsubsection{Example}\label{Yang-Mills extra-confluent}
We consider the presentation from Example~\ref{Yang-Mills side-confluent} of the Yang-Mills algebra over two generators. The vector space $\V{\otimes 2}\otimes\overline{R}\cap\overline{R}\otimes\V{\otimes 2}$ is reduced to $\{0\}$. Then, the extra-condition holds. We conclude that the presentation from Example~\ref{Yang-Mills side-confluent} is extra-confluent.

\section{Confluence algebras and reduction operators}\label{Confluence algebras and reduction operators}

\subsection{The contracting homotopy in small degree}\label{The contracting homotopy in small degree}

Through this section we fix an $N$-homogeneous algebra $A$. We assume that $A$ admits an extra-confluent presentation $\EV{X\mid R}$. This presentation is also fixed. We write $V=\K X$. 

The aim of this section is to make explicit our contracting homotopy in small homological degree. The formal construction will be done in Section~\ref{The left bound of a side-confluent presentation}.

We have to construct a family of $\K$-linear maps
\[h_{-1}:\K\F{} A,\ \text{and}\ h_n:A\otimes J_n\F{}A\otimes J_{n+1},\ \text{for}\ 0\leq n\leq 2,\]
satisfying the following relations:
\[\partial_1h_0+h_{-1}\varepsilon=\id{A}\ \text{and}\ \partial_{n+1}h_n+h_{n-1}\partial_n=\id{A\otimes J_n},\ \text{for}\ 0\leq n\leq 2.\]
By assumption, the set $X$ is finite. However, we will see that for the constructions of $h_{-1}$, $h_0$ and $h_1$ this hypothesis is not necessary.

From Proposition~\ref{Diamond lemma}, every element $f$ of $\tens{V}$ admits a unique normal form for $\EV{X\mid R}$. This normal form is denoted by $\widehat{f}$.

For every $w\in\EV{X}$, we define $[w]\in A\otimes V$ as follows:
\[[w]=\left\{\begin{split}
&0,\ \text{if}\ w\ \text{is the empty word},\\
&\overline{w'}\otimes x,\ \text{where}\ w'\in\EV{X}\ \text{and}\ x\in X\ \text{are such that}\ w=w'x.
\end{split}\right.\]
The map $[\ ]:\EV{X}\F{} A\otimes V$ is extended into a $\K$-linear map from $\tens{V}$ to $A\otimes V$. Let $w\in\EV{X}$ be a non empty word. For every $a\in A$, the action of $a$ on $[w]$ is given by $a.[w]=\left[fw\right]$, where $f\in\tens{V}$ is such that $a=\overline{f}$.

In small degree, the Koszul complex of $A$ is
\[A\otimes\left(V\otimes\vect{R}\cap\vect{R}\otimes V\right)\overset{\partial_3}{\F{}}A\otimes\overline{R}\overset{\partial_2}{\F{}}A\otimes V\overset{\partial_1}{\F{}}A\overset{\varepsilon}{\F{}}\K\F{}0,\]
where $\partial_1$ is defined by $\partial_1\left(1_A\otimes v\right)=\overline{v}$ for every $v\in V$, $\partial_2$ is defined by $\partial_2\left(1_A\otimes f\right)=\left[f\right]$ for every $f\in\overline{R}$ and $\partial_3$ is defined by $\partial_3\left(1_A\otimes g\right)=\sum\overline{v}\otimes f$ where $\sum vf$ is a decomposition of $g\in V\otimes\overline{R}\cap\overline{R}\otimes V$ in $V\otimes\overline{R}$. By definition of $\partial_3$, $\partial_3\left(1_A\otimes g\right)$ does not depend on the decomposition of $g$ in $V\otimes\overline{R}$.

\subsubsection{The constructions of $h_{-1}$ and $h_0$}

The maps $h_{-1}:\K\F{} A$ and $ h_0~:~A~\F{}~ A~\otimes ~V$ are defined by
\[h_{-1}\left(1_{\K}\right)=1_A\ \text{and}\ h_0(a)=\left[\widehat{f}\right],\ \text{where}\ f\in\tens{V}\ \text{is such that}\ \overline{f}=a.\]
We have $h_0\left(1_A\right)=0$ and $h_{-1}\varepsilon\left(1_A\right)=1_A$. If $a$ belongs to $A_m$ for $m\geq 1$, we have $\varepsilon\left(a\right)=0$ and $\partial_1h_0\left(a\right)=\overline{f}$. It follows that $\partial_1h_0+h_{-1}\varepsilon$ is equal to $\id{A}$.

\subsubsection{The construction of $h_{1}$}

Recall from Proposition~\ref{Diamond lemma} that the algebra $A$ admits as a basis the set $\left\{\overline{w},\ w\in\EV{X}\ \text{is a normal form}\right\}$. Thus, in order to define $h_1:A\otimes V\F{} A\otimes R$, it is sufficient to define $h_1\left(\overline{w}\otimes x\right)$ for every normal form $w\in\EV{X}$ and every $x\in X$. Moreover, $h_1$ has to satisfy the relation
\begin{equation} \label{equation of $h_1$}
\partial_2h_1\left(\overline{w}\otimes x\right)=\overline{w}\otimes x-h_0\left(\overline{wx}\right),
\end{equation}
for every normal form $w\in\EV{X}$ and every $x\in X$.

We define $h_1\left(\overline{w}\otimes x\right)$ by Noetherian induction on $wx$. Assume that $wx$ is a normal form. Then, let $h_1\left(\overline{w}\otimes x\right)=0$. We have:
\[\begin{split}
h_0\left(\overline{w}\overline{x}\right)&=[\widehat{wx}]\\
&=[wx]\\
&=\overline{w}\otimes x.
\end{split}\]
Thus, Relation~\ref{equation of $h_1$} holds. Assume that $wx$ is not a normal form and that $h_1\left(\overline{w'}\otimes x'\right)$ is defined and satisfies $(E_1)$ for every normal form $w'\in\EV{X}$ and every $x'\in X$ such that $w'x' < wx$. The word $wx$ can be written as a product $w_1w_2$ where $w_2\in X^{(N)}$ is not a normal form. The presentation $\EV{X\mid R}$ is reduced. Thus, there exists a unique $f\in R$ such that $f=w_2-\widehat{w_2}$. Let
\[h_1\left(\overline{w}\otimes x\right)=\overline{w_1}\otimes f+h_1\left([w_1\widehat{w_2}]\right).\]
We have:
\[\begin{split}
\partial_2h_1\left(\overline{w}\otimes x\right)&=[w_1f]+\partial_2h_1\left([w_1\widehat{w_2}]\right)\\
&=[w_1w_2]-[w_1\widehat{w_2}]+\partial_2h_1\left([w_1\widehat{w_2}]\right).
\end{split}\]
By induction hypothesis, $\partial_2h_1\left([w_1\widehat{w_2}]\right)$ is equal to $\left[w_1\widehat{w_2}\right]-\left[\widehat{w_1w_2}\right]$. Hence, we have:
\[\begin{split}
\partial_2h_1\left(\overline{w}\otimes x\right)&=\left[w_1w_2\right]-\left[\widehat{w_1w_2}\right]\\
&=\overline{w}\otimes x-\left[\widehat{wx}\right]\\
&=\overline{w}\otimes x-h_0\left(\overline{w}\overline{x}\right).
\end{split}\]
Thus, Relation~\ref{equation of $h_1$} holds.

\subsubsection{Remark}

We consider the $\K$-linear morphisms
\[\begin{split}
&F_1:A\otimes V\F{}\V{\otimes N},\ \overline{w_1w_2}\otimes x\longmapsto\overline{w_1}\otimes w_2x,\\
&F_1^1:A\otimes\V{\otimes N}\F{}A\otimes V,\ \overline{w_1}\otimes w_2x\longmapsto\overline{w_1w_2}\otimes x,\\
&F_2^1:A\otimes V\F{} A \otimes\V{\otimes N},\ \overline{w_1w_2}\otimes x\longmapsto\overline{w_1}\otimes\widehat{w_2x}.
\end{split}\]
The inductive definition of $h_1$ implies that $h_1\left(\overline{w}\otimes x\right)$ is equal to
\[\left(F_1-F_2^1\right)\left(\overline{w}\otimes x\right)+\left(F_1-F_2^1\right)\left(F_1^1F_2^1\left(\overline{w}\otimes x\right)\right)+\left(F_1-F_2^1\right)\left(\left(F_1^1F_2^1\right)^2\left(\overline{w}\otimes x\right)\right)+\cdots ,\]
where $\left(F_1-F_2^1\right)\left(\left(F_1^1F_2^1\right)^{2k}\left(\overline{w}\otimes x\right)\right)$ is vanishes for $k$ sufficiently large.\\\\

In order to define $h_2$ we need the following:

\subsubsection{Lemma}\label{Corollary extra-confluent and critical pair}

\textit{Let $A$ be an $N$-homogeneous algebra. Assume that $A$ admits an extra-confluent presentation $\EV{X\mid R}$. Let $w_1\in\EV{X}$, $w_2\in X^{\left(N-1\right)}$ and $x_1,x_2\in X$ such that}:
\begin{enumerate}
\item \textit{$w_1x_1$ and $x_1w_2$ are normal forms for $\EV{X\mid R}$.}
\item \textit{$w_2x_2$ is not a normal form for $\EV{X\mid R}$.}
\end{enumerate}
\textit{The word $w_1x_1w_2$ is a normal form for $\EV{X\mid R}$.}

\begin{proof}
Assume that $w_1x_1w_2$ is not a normal form. By hypothesis, $w_1x_1$ and $x_1w_2$ are normal forms. Thus, there exist a right divisor $u$ of $w_1$ and a left divisor $v$ of $w_2$ such that $ux_1v$ has length $N$ and is not a normal form. In particular, $ux_1w_2x_2$ is the source of a critical branching. From Proposition~~\ref{extra condition and critical pair}, the word $x_1w_2$ is not normal, which is a contradiction. Thus, Lemma~\ref{Corollary extra-confluent and critical pair} holds.
\end{proof}

\subsubsection{The construction of $h_2$}

Recall from Proposition~\ref{Diamond lemma} that the algebra $A$ admits as a basis the set $\left\{\overline{w},\ w\in\EV{X}\ \text{is a normal form}\right\}$. Thus, in order to define $h_2:A\otimes R\F{} A\otimes J_3$ it is sufficient to define $h_1\left(\overline{w}\otimes f\right)$ for every normal form $w\in\EV{X}$ and every $f\in R$. Moreover, $h_2$ has to satisfy the relation
\begin{equation} \label{equation of $h_2$}
\partial_3h_2\left(\overline{w}\otimes f\right)=\overline{w}\otimes f-h_1\partial_2\left(\overline{w}\otimes f\right),
\end{equation}
for every normal form $w\in\EV{X}$ and every $f\in R$.

We write $w=w_1x_1$, $f=w'-\widehat{w'}$ and $w'=w_2x_2$. We define $h_2\left(\overline{w}\otimes f\right)$ by Noetherian induction on $x_1w_2$. Assume that $x_1 w_2$ is a normal form. Let $h_2\left(\overline{w}\otimes f\right)=0$. We have:
\[\begin{split}
h_1\partial_2\left(\overline{w}\otimes f\right)&=h_1\left([wf]\right)\\
&=h_1\left([ww']\right)-h_1\left([w\widehat{w'}]\right)\\
&=h_1\left(\overline{ww_2}\otimes x'\right)-h_1\left([w\widehat{w'}]\right).
\end{split}\]
From Lemma~\ref{Corollary extra-confluent and critical pair}, $ww_2$ is a normal form. Thus, by construction of $h_1$, we have:
\[h_1\left(\overline{ww_2}\otimes x'\right)=\overline{w}\otimes f+ h_1\left([w\widehat{w'}]\right).\]
We conclude that $h_1\partial_2\left(\overline{w}\otimes f\right)$ is equal to $\overline{w}\otimes f$. Hence, Relation~\ref{equation of $h_2$} holds.

Assume that $h_2\left(\overline{u}\otimes g\right)$ is defined and that $(E_2)$ holds for every normal form $u\in\EV{X}$ and $g\in R$ such that $yv<x_1w_2$, where $y\in X$ and $v\in X^{(N-1)}$ are such that $u=u'y$ and $\lm{g}=vz$ for $u'\in\EV{X}$ and $z\in X$. We consider the two morphisms
\[S_1=S\otimes\id{V}\ \text{and}\ S_2=\id{V}\otimes S.\]
The presentation $\EV{X\mid R}$ is side-confluent. Thus, from Lemma~\ref{Lambda}, there exists an integer $k$ such that:
\[\EV{\id{\V{\otimes N+1}}-S_2,\id{\V{\otimes N+1}}-S_1}^k=\EV{\id{\V{\otimes N+1}}-S_1,\id{\V{\otimes N+1}}-S_2}^k.\]
We denote by $\Lambda$ this common morphism. The image of $\Lambda$ is included in $\ker\left(S_1\right)\cap\ker\left(S_2\right)$. The latter is equal to $\left(\overline{R}\otimes V\right)\cap\left(V\otimes\overline{R}\right)$. Recall that we have:
\[\EV{\id{\V{\otimes N+1}}-S_2,\id{\V{\otimes N+1}}-S_1}^k=\id{\V{\otimes N+1}}+\sum_{i=1}^{k-1}(-1)^i\left(\EV{S_1,S_2}^i+\EV{S_2,S_1}^i\right)+(-1)^{k}\EV{S_2,S_1}^k.\]
Thus, we have:
\[\Lambda=\left(\id{\V{\otimes N+1}}-S_2\right)+\left(\id{\V{\otimes N+1}}-S_2\right)\sum_{i=1}^{k-1}(-1)^i\ g_i\left(S_1,S_2\right),\]
where $g_i\left(S_1,S_2\right)$ denotes the product $S_1S_2S_1\cdots $ with $i$ factors. In particular, there exist $f_1,\cdots ,f_l\in R$, $x_1,\cdots ,x_l\in X$ and $\lambda_1,\cdots ,\lambda_l\in\K$ such that $x_iw_i<x_1w_2$ where $\lm{f_i}=w_iy_i$ and
\[\Lambda\left(xw'\right)=xf+\sum_{i=1}^l\lambda_ix_if_i.\]
Then, let
\[h_2\left(\overline{w}\otimes f\right)=\overline{w_1}\otimes\Lambda\left(xw'\right)-\sum_{i=1}^l\lambda_ih_2\left(\overline{w_1x_i}\otimes f_i\right).\]
We will show in Section~\ref{The left bound of a side-confluent presentation} that Relation~\ref{equation of $h_2$} holds.

\subsubsection{Remark}

We consider the $\K$-linear maps
\[\begin{split}
&F_2:A\otimes\V{\otimes N}\F{}\V{\otimes N+1},\ \overline{w_1x}\otimes w_2\longmapsto\overline{w_1}\otimes xw_2, \\
&F_1^2\colon A\otimes\V{\otimes N+1}\longrightarrow A\otimes\V{\otimes N},\ \overline{w_1}\otimes xw_2\longmapsto\overline{w_1x}\otimes w_2, \\
&F_2^2:A\otimes\V{\otimes N}\F{}A\otimes\V{\otimes N+1},\ \overline{w_1}x\otimes w_2\longmapsto\overline{w_1}\otimes xw_2-\Lambda\left(xw_2\right).
\end{split}\]
The inductive definition of $h_2$ implies that $h_2\left(\overline{w}\otimes f\right)$ is equal to
\[\left(F_2-F_2^2\right)\left(\overline{w}\otimes f\right)+\left(F_2-F_2^2\right)\left(F_2^1F_2^2\left(\overline{w}\otimes f\right)\right)+\left(F_2-F_2^2\right)\left(\left(F_2^1F_2^2\right)^2\left(\overline{w}\otimes f\right)\right)+\cdots ,\]
where $\left(F_2-F_2^2\right)\left(\left(F_2^1F_2^2\right)^{2k}\left(\overline{w}\otimes f\right)\right)$ is vanishes for $k$ sufficiently large.

\subsubsection{Example}

The construction of our contracting homotopy for the Koszul complex of the Yang-Mills algebra over two generators is done in Section~\ref{The Yang-Mills algebra over two generators}.

\subsection{Reduction operators and confluence algebras}\label{Reduction operators confluence algebras}

We fix a finite set $Y$, totally ordered by a relation $<$. For every $v\in\K Y\setminus\{0\}$, we denote by $\text{lm}(v)$ the greatest element of $Y$ occurring in the decomposition of $v$. We extend the order $<$ to a partial order on $\K Y$ in the following way: we have $v<w$ if either $v=0$ or if $\text{lm}(v)<\text{lm}(w)$.

In this section we recall some results from~\cite{MR1608711} about reductions operators and confluence algebras. 

\subsubsection{Reduction operators}\label{Reduction operators}

A linear projector $T$ of $\K Y$ is called a \textit{reduction operator relatively to} $Y$ if for every $y\in Y$, we have either $T(y)=y$ or $T(y)<y$. We denote by $\text{Red}\left(Y\right)$ the set of reduction operators relatively to $Y$.

\subsubsection{Lattice structure}\label{Lattice structure}

The set $\red{Y}$ admits a lattice structure. To define the order, recall from~\cite[Lemma 2.2]{MR1608711} that if $U,T\in\red{Y}$ are such that $\ker\left(U\right)$ is included in $\ker\left(T\right)$, then $\im{T}$ is included in $\im{U}$. Thus, the relation defined by $T\preceq U$ if $\ker(U)~\subset~\ker(T)$ is an order relation on $\red{Y}$.

We denote by $\M{L}(\K Y)$ the lattice of sub vector spaces of $\K Y$: the order is the inclusion, the lower bound is the intersection and the upper bound is the sum. To define the upper bound and the lower bound on $\red{Y}$, recall from~\cite[Theorem 2.3]{MR1608711}  that the map
\[\begin{split}
  \theta_{Y}\colon\red{Y}&\longrightarrow\M{L}(\K Y),\\
  T&\longmapsto\ker(T)
\end{split}
\] 
is a bijection. The lower bound $T_1\wedge T_2$ and the upper bound $T_1\vee T_2$ of two elements $T_1$ and $T_2$ of $\red{Y}$ are defined in the following way:
\[\left\{\begin{split}
&T_1\wedge T_2=\theta_Y^{-1}\left(\ker(T_1)+\ker(T_2)\right), \\ 
&T_1\vee T_2=\theta_Y^{-1}\left(\ker(T_1)\cap\ker(T_2)\right).
\end{split}
\right.\]

\subsubsection{Remark}\label{Minimum and maximum}

The lattice $\red{Y}$ admits $\id{\K Y}$ as maximum and $0_{\K Y}$ as minimum.

\subsubsection{Confluent pairs of reduction operators}\label{Confluent pairs of reduction operators}

A pair $P=\left(T_1,T_2\right)$ of reduction operators relatively to $Y$ is said to be \emph{confluent} if there exists an integer $k$ such that:
\[\EV{T_1,T_2}^k=\EV{T_2,T_1}^k.\]
We will see in Section~\ref{Reduction operators and side-confluent presentations} the link between this notion and the side-confluent presentations.

\subsubsection{Confluence algebras}\label{Confluence algebras}

Let $k$ be an integer. The \emph{confluence algebra of degree k} is the algebra presented by
\[\EV{s_1,s_2\mid\ s_i^2=s_i,\ \EV{s_1,s_2}^k=\EV{s_2,s_1}^k,\ i=1,2}.\]
This algebra is denoted by $\M{A}_k$. Let us consider the following elements of $\M{A}_k$:
\[\begin{split}
\sigma&=\EV{s_1,s_2}^k=\EV{s_2,s_1}^k,\\
\gamma_1&=(1-s_2)\sum_{i\in I}\EV{s_2,s_1}^i,\\
\gamma_2&=(1-s_1)\sum_{i\in I}\EV{s_1,s_2}^i,\\
\lambda&=1-\left(\sigma+\gamma_1+\gamma_2\right),
\end{split}\]
where $I$ is the set of odd integers between 1 and $k-1$. We easily check that we have the following relations:
\begin{subequations}
\begin{equation}\label{multiplication on the right}
\gamma_is_i=\gamma_i,\ \text{for}\ i=1,2,
\end{equation}
\begin{equation}\label{multiplication on the left}
s_i\gamma_i=s_i-\sigma,\ \text{for}\ i=1,2.
\end{equation}
\end{subequations}

\subsubsection{$P$-representations of confluence algebras}\label{$P$-representations of confluence algebras}

Let $P=\left(T_1,T_2\right)$ be a confluent pair of reduction operators relatively to $Y$. Let $k$ be an integer such that $\EV{T_1,T_2}^k=\EV{T_2,T_1}^k$. We consider the morphism of algebras
\[\begin{split}
  \varphi^P\colon\M{A}_{k}&\longrightarrow\text{End}\left(\K Y\right).\\
  s_i&\longmapsto T_i
\end{split}
\]
The morphism $\varphi^P$ is called the $P$-\emph{representation of} $\M{A}_k$. Recall from ~\cite{MR1608711} that:
\begin{subequations}
\begin{equation}\label{lower bound}
\varphi^P\left(\sigma\right)=T_1\wedge T_2,
\end{equation}
\begin{equation}\label{upper bound}
\varphi^P\left(1-\lambda\right)=T_1\vee T_2.
\end{equation}
\end{subequations}

\subsubsection{The left bound and the right bound}\label{The left bound}

Let $P=\left(T_1,T_2\right)$ be a confluent pair of reduction operators relatively to $Y$. By definition of $\lambda$ and from~\ref{$P$-representations of confluence algebras} we have:
\begin{equation} \label{Link between upper bound and lower bound}
T_1\vee T_2=T_1\wedge T_2+\varphi^P\left(\gamma_1\right)+\varphi^P\left(\gamma_2\right).
\end{equation}
The morphisms $\varphi^P\left(\gamma_1\right)$ and $\varphi^P\left(\gamma_2\right)$ are called the \emph{left bound} of $P$ and the \emph{right bound of P}, respectively. \\

We end this section with the following:

\subsubsection{Lemma}\label{Restriction}

\textit{Let $P=\left(T_1,T_2\right)$ be a confluent pair of reduction operators relatively to Y. Let $W$ be a sub vector space of $\K Y$. If $W$ is included in $\ker\left(T_i\right)$ for $i=1$ or 2, we have:}
\[{\varphi^P\left(\gamma_i\right)}_{\mid W}={T_1\vee T_2}_{\mid W}.\]

\begin{proof}
By definition, $\sigma$ and $\gamma_i$ factorize on the right by $s_i$. Hence, the restrictions of $\varphi^P\left(\sigma\right)$ and $\varphi^P\left(\gamma_i\right)$ to $W$ vanish. Thus, Lemma~\ref{Restriction} is a consequence of Relation~\ref{Link between upper bound and lower bound}.
\end{proof}

\subsection{Reduction operators and side-confluent presentations}\label{Reduction operators and side-confluent presentations}

Let $A$ be an $N$-homogeneous algebra. We suppose that $A$ admits a side-confluent presentation $\EV{X\mid R}$ where $X$ is a totally ordered finite set. For every integer $m$, the set $X^{(m)}$ is finite and totally ordered for the lexicographic order induced by the order over $X$. We write $V=\K X$.

\subsubsection{Normal forms and the Koszul complex}\label{Reduction operators and normal forms}

In Lemma~\ref{reduction operators and normal forms} we will link together the Koszul complex of $A$ and the reduction operators. In this way, recall from Proposition~\ref{Diamond lemma} that every element $f\in\tens{V}$ admits a unique normal form for $\EV{X\mid R}$, denoted by $\widehat{f}$. Let
\[\begin{split}
  \phi\colon\tens{V}&\longrightarrow\tens{V},\\
  f&\longmapsto\widehat{f}
\end{split}
\]
Recall from Proposition~\ref{Diamond lemma} that for every $f\in\tens{V}$, we have $f\in I(R)$ if and only if $\widehat{f}=0$. Hence, $\phi$ induces a $\K$-linear isomorphism $\overline{\phi}$ from $A$ to $\im{\phi}$. In particular, for every integer $n$, the morphism $\phi_n=\overline{\phi}\otimes\id{\V{\otimes l_N(n)}}$ is a $\K$-linear isomorphism from $A\otimes J_n$ to $\im{\phi}\otimes J_n$. Thus, the Koszul complex $\left(K_\bullet,\partial\right)$ of $A$ is isomorphic to the complex of vector spaces $\left(K'_\bullet,\partial'\right)$
\[\cdots\overset{\partial'_{n+1}}{\F{}}\im{\phi}\otimes J_n\overset{\partial'_n}{\F{}}\im{\phi}\otimes J_{n-1}\F{}\cdots\overset{\partial'_2}{\F{}}\im{\phi}\otimes J_1\overset{\partial'_1}{\F{}}\im{\phi}\overset{\varepsilon'}{\F{}}\K\F{}0,\]
where $\partial_n'$ is equal to $\phi_{n-1}\circ\partial_n\circ\phi_n^{-1}$.

\subsubsection{Definition}\label{Reduction operators and the Koszul complex}

The complex $\left(K'_\bullet,\partial'\right)$ is the \emph{normalised Koszul complex of A.}

\subsubsection{Lemma}\label{reduction operators and normal forms}

\begin{enumerate}
\item\label{normalisation operator and reduction operators} \textit{For every integer $m$, the restriction of $\phi$ to $\V{\otimes m}$ is a reduction operator relatively to $X^{(m)}$ and its kernel is equal to $I(R)_m$.}
\item\label{reduction operators and differentials}
\textit{Let $n$ be an integer such that $n\geq 1$. The morphism $\partial'_n$ is the restriction to $\emph{im}\left(\phi\right)\otimes J_n$ of the morphism $\varphi_n:\bigoplus_{m\geq l_N(n)}\V{\otimes m}\F{}\emph{T}\left(V\right)$ defined by}
\[{\varphi_n}_{\mid\V{\otimes m}}=\phi_{\mid\V{\otimes m-l_N(n-1)}}\otimes\id{\V{\otimes l_N(n-1)}}.
\]
\end{enumerate}
 
\begin{proof}
Let us show the first point. The presentation $\EV{X\mid R}$ is $N$-homogeneous. Thus, for every $w\in X^{(m)}$, $\phi(w)$ belongs to $\V{\otimes m}$. In particular, the restriction of $\phi$ to $\V{\otimes m}$ is an endomorphism of $\V{\otimes m}$. Let $w\in X^{(m)}$. If $w$ is a normal form, then $\phi(w)$ is equal to $w$. In particular, $\phi_{\mid\V{\otimes m}}$ is a projector. If $w$ is not a normal form, then $\phi(w)=\widehat{w}$ is strictly smaller than $w$. Thus, $\phi_{\mid\V{\otimes m}}$ is a reduction operator relatively to $X^{(m)}$. Moreover, $\widehat{f}$ is equal to $0$ if and only if $f$ belongs to $I(R)$. Thus, the kernel of $\phi_{\mid\V{\otimes m}}$ is equal to $I(R)_m$.

Let us show the second point. Recall from~\ref{The construction of the Koszul complex} that the differential $\partial_n~:~A~\otimes ~J_n~\F{}~ A\otimes~ J_{n-1}$ of the Koszul complex of $A$ is the restriction to $A\otimes J_n$ of the $A$-linear map defined by:
\[
\begin{split}
 A\otimes\V{\otimes l_N(n)}&\F{} A\otimes\V{\otimes l_N(n-1)},\\
1_A\otimes w&\longmapsto \overline{w_1}\otimes w_2
\end{split}
\]
where $w_1\in X^{\left(l_N(n)-l_N(n-1)\right)}$ and $w_2\in X^{\left(l_N(n-1)\right)}$ are such that $w=w_1w_2$. Thus, the map $\partial'_n$ is the restriction of the morphism which maps a word $w$ of length $m\geq l_N(n)$ to $\widehat{w_1}w_2$, where $w_1~\in~ X^{\left(m-l_N(n-1)\right)}$ and $w_2\in X^{\left(l_N(n-1)\right)}$ are such that $w=w_1w_2$. The latter is equal to $\phi_{\mid\V{\otimes m-l_N(n-1)}}\otimes\id{\V{\otimes l_N(n-1)}}$.
\end{proof}

\subsubsection{Lattice properties}\label{Lattice properties}

Let $S\in\text{End}\left(\V{\otimes N}\right)$ be the operator of $\EV{X\mid R}$:
\[S\left(w\right)=\left\{
\begin{split}
&\lm{f}-f,\ \text{if there exists}\ f\in R\ \text{such that}\ w=\lm{f}, \\
&w,\ \text{otherwise}.
\end{split}
\right.\]
The properties of $S$ described in~\ref{The operator of a presentation} imply that $S$ is equal to $\theta_{X^{(N)}}^{-1}\left(\vect{R}\right)$. For every integers $m$ and $i$ such that $m\geq N$ and $0\leq i\leq m-N$, we consider the following reduction operator relatively to $X^{(m)}$:
\[S^{(m)}_i=\id{\V{\otimes i}}\otimes S\otimes\id{\V{\otimes m-N-i}}.\]
The kernel of $S^{(m)}_i$ is equal to $\V{\otimes i}\otimes\overline{R}\otimes\V{\otimes m-N-i}$.

The presentation $\EV{X\mid R}$ is side-confluent. Hence, the pair $\left(S_i^{\left(2N-1\right)},S_j^{\left(2N-1\right)}\right)$ is confluent for every integers $i$ and $j$ such that $0~\leq ~i~,~j~\leq ~N~-~1$. We deduce from \cite[Section 3]{MR1832913} and \cite[Theorem 2.12]{MR1608711} that for every integer $m$ such that $m\geq N$, the sub-lattice of $\red{X^{(m)}}$ spanned by $S^{(m)}_0,\cdots ,S^{(m)}_{m-N}$ is \emph{confluent} (that is, the elements of this lattice are pairwise confluent) and \emph{distributive} (that is, for every $S,T,U$ belonging to this lattice, we have $\left(S\wedge T\right)\vee U=\left(S\vee U\right)\wedge\left(T\vee U\right)$).

\section{The left bound of a side-confluent presentation}\label{The left bound of a side-confluent presentation}

Through this section we fix an $N$-homogeneous algebra $A$. We assume that $A$ admits an $N$-homogeneous presentation $\EV{X\mid R}$ where $X$ is a totally ordered finite set. This presentation is also fixed. We write $V=\K X$. We consider the notations of~\ref{Lattice properties}.

\subsection{Reduction pairs associated with a presentation}\label{Reduction pairs associated to a presentation}

For every integers $n$ and $m$ such that $m\geq l_N(n)$, we consider the following reduction operators relatively to $X^{(m)}$:
\[{F_1^{n,m}}=\theta_{X^{(m)}}^{-1}\left(I(R)_{m-l_N(n)}\otimes\V{\otimes l_N(n)}\right),\]
\[F_2^{n,m}=\left\{\begin{split}
&\id{V{\otimes m}},\ \text{if}\  m <l_N(n+1),\\
&\theta_{X^{(m)}}^{-1}\left(\V{\otimes m-l_N(n+1)}\otimes J_{n+1}\right),\ \text{otherwise}.
\end{split}
\right.\]
We denote by $P_{n,m}$ the pair $\left(F_1^{n,m},F_2^{n,m}\right)$.

\subsubsection{Definition}

The pair $P_{n,m}$ is the \emph{reduction pair of bi-degree (n,m) associated with} $\EV{X\mid R}$.

\subsubsection{Lemma}\label{parity} 

\textit{Let $n$ and $m$ be two integers such that $n\geq 1$ and $l_N(n)\leq m< l_N(n+1)$. Then, $m-l_N(n-1)$ is smaller than $N-1$ and $F_1^{n-1,m}$ is equal to $\idd{\V{\otimes m}}$.}

\begin{proof}
First, we show that $m-l_N(n-1)$ is smaller than $N-1$. Assume that $m$ is a multiple of $N$: $m=kN$. In this case, the hypothesis $l_N(n)\leq m< l_N(n+1)$ implies that $n$ is equal to $2k$. Thus, $l_N(n-1)$ is equal to $(k-1)N+1$. That implies that $m-l_N(n-1)$ is equal to $N-1$. Assume that $m$ is not a multiple of $N$: $m=kN+r$ with $1\leq r\leq N-1$. In this case, the hypothesis $l_N(n)\leq m< l_N(n+1)$ implies that $n$ is equal to $2k+1$. Thus, $m-l_N(n-1)=m-kN$ is smaller than $N-1$. 

Let us show that $F_1^{n-1,m}$ is equal to $\id{\V{\otimes m}}$. The first part of the lemma implies that $I(R)_{m-l_N(n-1)}$ is equal to $\{0\}$. Thus, the kernel of $F_1^{n-1,m}$ is equal to $\{0\}$, that is, $F_1^{n-1,m}$ is equal to $\id{\V{\otimes m}}$.
\end{proof}

\subsubsection{Lemma}\label{appartence au treillis}
\begin{enumerate}
\item \textit{Let $n$ and let $m$ be two integers such that $m\geq l_N(n+2)$. We have:}
\[F^{n,m}_1=S^{(m)}_0\wedge\cdots\wedge S^{(m)}_{m-l_N(n+2)}.\]
\item \textit{Let $n$ and $m$ be two integers such that $n\geq 1$ and $m\geq l_N(n+1)$. We have:}
\[F^{n,m}_2=S^{(m)}_{m-l_N(n+1)}\vee\cdots\vee S^{(m)}_{m-N}.\]
\end{enumerate}

\begin{proof}
By definition of $\wedge$, we have:
\[\begin{split}
\ker\left(S^{(m)}_0\wedge\cdots\wedge S^{(m)}_{m-l_N(n+2)}\right)&=\sum_{i=0}^{m-l_N(n+2)}\ker\left(S^{(m)}_i\right) \\
&=\sum_{i=0}^{m-l_N(n+2)}\V{\otimes i}\otimes\overline{R}\otimes\V{\otimes m-N-i}\\
&=\left(\sum_{i=0}^{m-l_N(n+2)} \V{\otimes i}\otimes\overline{R}\otimes\V{\otimes m-l_N(n)-N-i}\right)\otimes\V{\otimes l_N(n)}\\
&=\left(\sum_{i=0}^{m-l_N(n)-N} \V{\otimes i}\otimes\overline{R}\otimes\V{\otimes m-l_N(n)-N-i}\right)\otimes\V{\otimes l_N(n)}\\
&=I(E)_{m-l_N(n)}\otimes\V{\otimes l_N(n)}.
\end{split}\]
By definition of $\vee$, we have:
\[\begin{split}
\ker\left(S^{(m)}_{m-l_N(n+1)}\vee\cdots\vee S^{(m)}_{m-N}\right)&=\bigcap_{i=m-l_N(n+1)}^{m-N} \ker\left(S^{(m)}_i\right)\\
&=\bigcap_{i=m-l_N(n+1)}^{m-N}\V{\otimes i}\otimes\overline{R}\otimes\V{\otimes m-N-i}\\
&=\V{\otimes m-l_N(n+1)}\otimes\left(\bigcap_{i=0}^{l_N(n+1)-N}\V{\otimes i}\otimes\overline{R}\otimes\V{\otimes l_N(n+1)-N-i}\right)\\
&=\V{\otimes m-l_N(n+1)}\otimes J_{n+1}.
\end{split}\]
The map $\theta_{X^{(m)}}$ being a bijection, the two relations hold.
\end{proof}
 
\subsubsection{Theorem}\label{Confluence}

\textit{Let $A$ be an $N$-homogeneous algebra admitting a side-confluent presentation $\EV{X\mid R}$, where $X$ is a finite set. The reduction pairs associated with $\EV{X\mid R}$ are confluent.}

\begin{proof}
Let $n$ and $m$ be two integers such that $l_N(n)\leq m$. We have to show that the reduction pair of bi-degree $(n,m)$ associated with $\EV{X\mid R}$ is confluent.

\paragraph{Step 1.} Assume that $n=0$. We have $P_{0,0}=\left(\id{\K},\id{\K}\right)$. Thus, the pair $P_{0,0}$ is confluent. Let $m$ be an integer such that $m\geq 1$. The kernel of $F^{0,m}_2$ is equal to $\V{\otimes m-1}\otimes J_1=\V{\otimes m}$. Thus, $F^{0,m}_2$ is equal to $0_{\V{\otimes m}}$. In particular, the operators $F_1^{0,m}$ and $F_2^{0,m}$ commute, that is, they satisfy the relation $\EV{F_1^{0,m},F_2^{0,m}}^2=\EV{F_2^{0,m},F_1^{0,m}}^2$. Hence, the pair $P_{0,m}$ is confluent for every integer $m$.

\paragraph{Step 2.} Assume that $n\geq 1$ and $l_n(n)\leq m< l_N(n+1)$. The pair $P_{n,m}$ is equal to $\left(F_1^{n,m},\id{\V{\otimes m}}\right)$. Thus, the operators $F_1^{n,m}$ and $F_2^{n,m}$ commute. We conclude that the pairs $P_{n,m}$ such that $n\geq 1$ and $l_n(n)\leq m< l_N(n+1)$ are confluent.

\paragraph{Step 3.} Assume that $n\geq 1$ and $l_N(n+1)\leq m< l_N(n+2)$. From Lemma~\ref{parity}, the morphism $F^{n,m}_1$ is equal to $\id{\V{\otimes m}}$. In particular, the operators $F^{n,m}_1$ and $F_2^{n,m}$ commute. Thus, the pairs $P_{n,m}$ such that $n\geq 1$ and $l_N(n+1)\leq m< l_N(n+2)$ are confluent.

\paragraph{Step 4.} Assume that $n\geq 1$ and $m\geq l_N(n+2)$. Lemma~\ref{appartence au treillis} implies that $F^{n,m}_1$ and $F^{n,m}_2$ belong to the lattice generated by $S^{(m)}_i$, for $0\leq i\leq m-N$. From~\ref{Lattice properties} that the latter is confluent. Hence, the pairs $P_{n,m}$ such that $n\geq 1$ and $m\geq l_N(n+2)$ are confluent.
\end{proof}

\subsection{Construction}\label{Construction}

Through this section, we assume that the presentation $\EV{X\mid R}$ of $A$ is side-confluent. From Proposition~\ref{Diamond lemma}, every element $f$ of $\tens{V}$ admits a unique normal for $\EV{X\mid R}$. This normal form is denoted by $\widehat{f}$. We denote by $\phi$ the endomorphism of $\tens{V}$ which maps an element to its unique normal form. We consider the notations of Section~\ref{Reduction pairs associated to a presentation}.

\subsubsection{Lemma}\label{expression $F_1$}

\textit{For every integers $n$ and $m$ such that $m\geq l_N(n)$, the operator $F_1^{n,m}$ is equal to $\phi_{\mid\V{\otimes m-l_N(n)}}\otimes\idd{\V{\otimes l_N(n)}}$.}

\begin{proof}
From Point~\ref{normalisation operator and reduction operators} of Lemma~\ref{reduction operators and normal forms}, the operator $\phi_{\mid\V{\otimes m-l_N(n)}}\otimes\id{\V{\otimes l_N(n)}}$ is a reduction operator relatively to $X^{(m)}$ and its kernel is equal to $I(R)_{m-l_n(n)}\otimes\V{\otimes l_n(n)}$. The map $\theta_{X^{(m)}}$ being a bijection, Lemma~\ref{expression $F_1$} holds.
\end{proof}

\subsubsection{Lemma}\label{definition of the family}

\textit{Let $n$ be an integer. Let $h'_n:\bigoplus_{m\geq l_N(n)}\V{\otimes m}\F{}\emph{T}\left(V\right)$ be the $\K$-linear map defined by}
\[{h'_n}_{\mid\V{\otimes m}}=\varphi^{P_{n,m}}\left(\gamma_1\right),\]
\textit{where $\varphi^{P_{n,m}}\left(\gamma_1\right)$ is the left bound of $P_{n,m}$. The image of $h'_n$ is included in $\imm{\phi}~\otimes~ J_{n+1}$.}

\begin{proof}
Let $m$ be an integer such that $m\geq l_N(n)$. By definition of the left bound, there exists an endomorphism $H$ of $\V{\otimes m}$ such that
\[\varphi^{P_{n,m}}\left(\gamma_1\right)=\left(\id{\V{\otimes m}}-F^{n,m}_2\right)F^{n,m}_1H.\]

The image of $F^{n,m}_1=\phi_{\mid\V{\otimes m-l_N(n)}}\otimes\V{\otimes l_N(n)}$ is equal to the vector space spanned by the elements with shape $w_1w_2$ where $w_1\in X^{\left(m-l_N(n)\right)}$ is a normal form and $w_2\in X^{\left(l_N(n)\right)}$.

Let
\[G=\theta_{X^{\left(l_N(n+1)\right)}}^{-1}\left(J_{n+1}\right).\]
We have $F^{n,m}_2=\id{\V{\otimes m-l_N(n+1)}}\otimes G$. The latter implies that 
\[\left(\id{\V{\otimes m}}-F^{n,m}_2\right)=\id{\V{\otimes m-l_N(n+1)}}\otimes\left(\id{\V{\otimes l_N(n+1)}}-G\right).\]

We conclude that the image of $\varphi^{P_{n,m}}\left(\gamma_1\right)$ is included in the vector space spanned by elements with shape $wf$ where $w\in X^{\left(m-l_N(n+1)\right)}$ is a normal form and $f\in J_{n+1}$. This vector space is equal to $\im{\phi_{\mid\V{\otimes m-l_N(n+1)}}}\otimes J_{n+1}$. 
\end{proof}

\subsubsection{Definition}\label{Definition of the left bound}

For every integer $n$, let
\[h_n=\phi_{n+1}^{-1}\circ h'_n\circ\phi_n:A\otimes J_n\F{}A\otimes J_{n+1},\]
where $\phi_n$ is the $\K$-linear isomorphism between $A\otimes J_n$ and $\im{\phi}\otimes J_n$ defined in~\ref{Reduction operators and normal forms}. The family $\left(h_n\right)_n$ is the \emph{left bound of} $\EV{X\mid R}$.

\subsubsection{Reduction relations}\label{Reduction relations}

Let $n$ and $m$ be two integers such that $m\geq l_N(n)$. Then, we denote by $K_n^{(m)}~=~\im{\phi_{\mid\V{\otimes m-l_N(n)}}}\otimes~ J_n$. In particular, we have:
\[\im{\phi}\otimes J_n=\bigoplus_{m\geq l_N(n)}K_n^{(m)}.\]
We say that the presentation $\EV{X\mid R}$ satisfy the \emph{reduction relations} if for every integers $n$ and $m$ such that $m\geq l_n(n)$, the following equality holds:
\[\left(r_{n,m}\right)\ \ \ {F_1^{n,m}\wedge F_2^{n,m}}_{\mid K^{(m)}_n}={F_1^{n-1,m}\vee F_2^{n-1,m}}_{\mid K^{(m)}_n}.\]

\subsubsection{Proposition}\label{relation r}

\textit{Let $A$ be an $N$-homogeneous algebra. Assume that $A$ admits a side-confluent presentation $\EV{X\mid R}$ where $X$ is a finite set. The left bound of $\EV{X\mid R}$ is a contracting homotopy for the Koszul complex of $A$ if and only if $\EV{X\mid R}$ satisfy the reduction relations.}

\begin{proof}
The left bound of $\EV{X\mid R}$ is a contracting homotopy for the Koszul complex of $A$ if and only if the family $\left(h'_n:\im{\phi}\otimes J_n\F{}\im{\phi}\otimes J_{n+1}\right)_n$ defined in Lemma~\ref{definition of the family} is a contracting homotopy for the normalised Koszul complex of $A$.

From Point~\ref{reduction operators and differentials} of Lemma~\ref{reduction operators and normal forms}, the restriction of $F_1^{n-1,m}=\phi_{\mid\V{\otimes m-l_N(n-1)}}\otimes\id{\V{\otimes l_n(n-1)}}$ to $K^{(m)}_n$ is equal to the restriction of $\partial'_n$ to $K^{(m)}_n$. Thus, the family $\left(h'_n\right)_n$ is a contracting homotopy for $\left(K'_\bullet,\partial'\right)$ if and only if for every $n$ and $m$ such that $n\geq 1$ and $m\geq l_N(n)$, the following relation holds:
\[\left(\varphi^{P_{n,m}}\left(s_1\right)\varphi^{P_{n,m}}\left(\gamma_1\right)+\varphi^{P_{n-1,m}}\left(\gamma_1\right)\varphi^{P_{n-1,m}}(s_1)\right)_{\mid K^{(m)}_n}=\id{K_n^{(m)}}.\]
From Relation~\ref{multiplication on the left} (see page~\pageref{multiplication on the left}) and Relation~\ref{lower bound} (see page~\pageref{lower bound}), we have:
\[\begin{split}
{\varphi^{P_{n,m}}\left(s_1\right)\varphi^{P_{n,m}}\left(\gamma_1\right)}&=F_1^{n,m}- \varphi^{P_{n,m}}(\sigma)\\
&={F_1^{n,m}-F_1^{n,m}\wedge F_2^{n,m}}.
\end{split}\]
The image of $F_1^{n,m}=\phi_{\mid\V{\otimes m-l_N(n)}}\otimes\id{\V{\otimes l_N(n)}}$ is equal to $\im{\phi_{\mid\V{\otimes m-l_N(n)}}}\otimes\V{\otimes l_N(n)}$. Thus, $K^{(m)}_n$ is included in $\im{F_1^{n,m}}$. In particular, the restriction of $F_1^{n,m}$ to $K^{(m)}_n$ is the identity map. We deduce that the left bound family of $\EV{X\mid R}$ is a contracting homotopy for the Koszul complex of $A$ if and only if the following relation holds:
\[\left(\varphi^{P_{n-1,m}}\left(\gamma_1\right)\varphi^{P_{n-1,m}}(s_1)\right)_{\mid K^{(m)}_n}={F_1^{n,m}\wedge F_2^{n,m}}_{\mid K^{(m)}_n}.\]
From Relation~\ref{multiplication on the right} (see page~\pageref{multiplication on the right}), $\varphi^{P_{n-1,m}}\left(\gamma_1\right)\varphi^{P_{n-1,m}}(s_1)$ is equal to $\varphi^{P_{n-1,m}}\left(\gamma_1\right)$. Thus, it is sufficient to show:

\begin{equation} \label{left expression}
\varphi^{P_{n-1,m}}\left(\gamma_1\right)_{\mid K^{(m)}_n}={F_1^{n-1,m}\vee F_2^{n-1,m}}_{\mid K^{(m)}_n}.
\end{equation}\noindent
By construction, $K^{(m)}_n$ is included in $\ker\left(F_2^{n-1,m}\right)$. Hence, Relation~\ref{left expression} is a consequence of Lemma~\ref{Restriction}.

\end{proof}

The following lemma will be used in the proof of Theorem~\ref{Theorem}:

\subsubsection{Lemma}\label{commute}

\textit{Let $n$ and $m$ be two integers such that $n\geq 1$ and $l_N(n)\leq m< l_N(n+1)$. The operators $F_1^{n,m}$ and $F_1^{n-1,m}\vee F_2^{n-1,m}$ commute.}

\begin{proof}
The pair $P_{n,m}$ being confluent, we deduce from Relation~\ref{upper bound} (see page~\pageref{upper bound}) that $F^{n-1,m}_1\vee F^{n-1,m}_2$ is polynomial in $F^{n-1,m}_1$ and $F^{n-1,m}_2$. Hence, it is sufficient to show that $F^{n,m}_1$ commutes with $F^{n-1,m}_1$ and $F^{n-1,m}_2$.

Let 
\[G=\theta_{X^{\left(l_N(n)\right)}}^{-1}\left(J_{n}\right).\]
We have $F^{n-1,m}_2=\id{\V{\otimes m-l_N(n)}}\otimes G$. Thus, $F_1^{n,m}=\phi_{\mid\V{\otimes m-l_N(n)}}\otimes\id{\V{\otimes l_N(n)}}$ commutes with $F^{n,m}_2$. Moreover, the morphism $F^{n,m}_1$ (respectively $F^{n-1,m}_1$) maps a word $w$ of length $m$ to $\widehat{w_1}w_2$ (respectively $\widehat{w'_1}w'_2$), where $w_1\in X^{\left(m-l_N(n)\right)}$ and $w_2\in X^{\left(l_N(n)\right)}$ (respectively $w'_1\in X^{\left(m-l_N(n-1)\right)}$ and $w'_2\in X^{\left(l_N(n-1)\right)}$) are such that $w=w_1w_2$ (respectively $w=w'_1w'_2$). Thus, the two compositions $F_1^{n,m}F_1^{n-1,m}$ and $F_1^{n-1,m}F_1^{n,m}$ are equal to $F_1^{n-1,m}$.
\end{proof}

\subsection{Extra-confluent presentations and reduction relations}\label{Extra-confluent presentattions and reduction relations}

Through this section we assume that the presentation $\EV{X\mid R}$ is extra-confluent. Our aim is to show that $\EV{X\mid R}$ satisfy the reduction relations. We consider the notations of Section~\ref{Reduction pairs associated to a presentation}.

\subsubsection{Lemma}\label{interpretation $(*)$}
\textit{Let $m$, $r$ and $k$ be three integers such that $m\geq N+2$, $2\leq k\leq N-1$ and  $r+k\leq m-N$. Then, we have:}
\begin{enumerate}
\item\label{point 1} $S^{(m)}_r\vee S^{(m)}_{r+k}=S^{(m)}_r\vee\cdots\vee S^{(m)}_{r+k}$,
\item\label{point 2} $\left(S^{(m)}_r\wedge\cdots\wedge S^{(m)}_{r+k-1}\right)\vee S^{(m)}_{r+k}=S^{(m)}_{r+k-1}\vee S^{(m)}_{r+k}$.
\end{enumerate}

\begin{proof}
Let us prove the point~\ref{point 1}. The extra-condition implies the following inclusion:
\[\left(\V{\otimes r+k}\otimes\overline{R}\otimes\V{\otimes m-N-r-k}\right)\cap\left(\V{\otimes r}\otimes\overline{R}\otimes\V{\otimes m-N-r}\right)\ \subset\ \V{\otimes r+k-1}\otimes\overline{R}\otimes\V{\otimes m-N-r-k+1}.\]
Applying the bijection $\theta_{X^{(m)}}^{-1}$, we have:
\[S^{(m)}_{r+k-1}\preceq S^{(m)}_r\vee S^{(m)}_{r+k}.\]
By definition of the upper bound, we deduce that $S^{(m)}_r\vee S^{(m)}_{r+k-1}\vee S^{(m)}_{r+k}$ is equal to $S^{(m)}_r\vee S^{(m)}_{r+k}$. By induction on $k$, we obtain the first relation.

Let us prove the point 2. Recall from~\ref{Lattice properties} that the lattice spanned by $S^{(m)}_0,\cdots ,S^{(m)}_{m-N}$ is distributive. Thus, the left hand side of 2 is equal to $\left(S^{(m)}_r\vee S^{(m)}_{r+k}\right)\wedge\cdots  \wedge\left(S^{(m)}_{r+k-1}\vee S^{(m)}_{r+k}\right)$. By the first point, for every $0\leq i\leq n-2$, $S^{(m)}_{r+i}\vee S^{(m)}_{r+k}$ is equal to $S^{(m)}_{r+i}\vee\cdots\vee S^{(m)}_{r+k}$, so it is greater than $S^{(m)}_{r+k-1}\vee S^{(m)}_{k+r}$. By definition of the lower bound, the second relation holds.
\end{proof}

\subsubsection{Lemma}\label{lemme intermediaire 1}

\textit{Let $n$ and $m$ be two integers such that $n\geq 2$ and $l_N(n+1)\leq m<l_N(n+2)$. We have:}
\begin{equation} \label{equation 1}
\left(S^{(m)}_0\wedge\cdots\wedge S^{(m)}_{m-l_N(n+1)}\right)\vee S^{(m)}_{m-l_N(n)}=S^{(m)}_{m-l_N(n+1)}\vee\cdots\vee S^{(m)}_{m-l_N(n)}.
\end{equation}

\begin{proof}
From Lemma~\ref{parity}, the hypothesis $l_N(n+1)\leq m<l_N(n+2)$ implies that $m-l_N(n)$ is smaller than $N-1$.

Assume that $m$ is a multiple of $N$. The hypothesis $l_N(n+1)\leq m<l_N(n+2)$ implies that $m$ is equal to $l_N(n+1)$. Thus, the left hand side of~\ref{equation 1} is equal to $S^{(m)}_0\vee S^{(m)}_{m-l_N(n)}$ and the right hand side of~\ref{equation 1} is equal to $S^{(m)}_0\vee\cdots\vee S^{(m)}_{m-l_N(n)}$. Hence, Relation~\ref{equation 1} is a consequence of Lemma~\ref{interpretation $(*)$} point~\ref{point 1}.

Assume that $m$ is not a multiple of $N$. The hypothesis $l_N(n+1)\leq m<l_N(n+2)$ implies that $n$ is even. Hence, the left hand side of~\ref{equation 1} is equal to $\left(S^{(m)}_0\wedge\cdots\wedge S^{(m)}_{m-l_N(n)-1}\right)\vee S^{(m)}_{m-l_N(n)}$ and the right hand side of~\ref{equation 1} is equal to $S^{(m)}_{m-l_N(n)-1}\vee S^{(m)}_{m-l_n(n)}$. If $n$ is equal to 2 and $m$ is equal to $N+1$, these two morphisms are equal to $S_0^{(N+1)}\vee S_1^{(N+1)}$. If the couple $(n,m)$ is different from $(2,N+1)$, Relation~\ref{equation 1} is a consequence of Lemma~\ref{interpretation $(*)$} point~\ref{point 2}.

\end{proof}

\subsubsection{Lemma}\label{lemme intermediaire 2}

\textit{Let $n$ and $m$ be two integers such that $n\geq 2$ and $m\geq l_N(n+2)$. Letting}
\[T_{n,m}=S^{(m)}_{m-l_N(n+2)+1}\wedge\cdots \wedge S^{(m)}_{m-l_N(n+1)},\]
\textit{we have:}
\[T_{n,m}\vee F_2^{n-1,m}=F_2^{n,m}.\]

\begin{proof}
From Lemma~\ref{appartence au treillis}, we have
\[\begin{split}
F_2^{n-1,m}&=S^{(m)}_{m-l_N(n)}\vee\cdots\vee S^{(m)}_{m-N} ,\ \text{and}\\
F_2^{n,m}&=S^{(m)}_{m-l_N(n+1)}\vee\cdots\vee S^{(m)}_{m-N}.
\end{split}\]
The law $\vee$ being associative, it is sufficient to show:
\begin{equation} \label{equation 2}
T_{n,m}\vee S^{(m)}_{m-l_N(n)}=S^{(m)}_{m-l_N(n+1)}\vee\cdots\vee S^{(m)}_{m-l_N(n)}.
\end{equation}

Assume that $n$ is odd. We have $l_N(n+2)=l_N(n+1)+1$. Hence, the left hand side of~\ref{equation 2} is equal to $S^{(m)}_{m-l_N(n+1)}\vee S^{(m)}_{m-l_N(n)}$. Moreover, $l_N(n+1)-l_N(n)$ is equal to $N-1$. Thus, Relation~\ref{equation 2} is a consequence of Lemma~\ref{interpretation $(*)$} point~\ref{point 1}. 

Assume that $n$ is even. We have $l_N(n+1)=l_N(n)+1$. Hence, the left hand side of~\ref{equation 2} is equal to $\left(S^{(m)}_{m-l_N(n+2)+1}\wedge\cdots \wedge S^{(m)}_{m-l_N(n)-1)}\right)\vee S^{(m)}_{m-l_N(n)}$ and the right hand side of~\ref{equation 2} is equal to $S^{(m)}_{m-l_n(n)-1}~\vee ~S^{(m)}_{m-l_N(n)}$. Moreover, $l_N(n+2)-1-l_N(n)$ is equal to $N-1$. Thus, Relation~\ref{equation 2} is a consequence of Lemma~\ref{interpretation $(*)$} point~\ref{point 2}.
\end{proof}

\subsubsection{Proposition}\label{utilisation $(*)$}

\textit{Let $A$ be an $N$-homogeneous algebra. Assume that $A$ admits an extra-confluent presentation $\EV{X\mid R}$. For every integers $n$ and $m$ such that $n\geq 1$ and $m\geq l_N(n+1)$, we have:}
\[F^{n,m}_1\wedge\left(F^{n-1,m}_1\vee F^{n-1,m}_2\right)=F_1^{n,m}\wedge F_2^{n,m}.\]

\begin{proof}
For every integers $n$ and $m$ such that $n\geq 1$ and $m\geq l_N(n+1)$, let
\[\begin{split}
L_{n,m}&=F^{n,m}_1\wedge\left(F^{n-1,m}_1\vee F^{n-1,m}_2\right), \\
R_{n,m}&=F_1^{n,m}\wedge F_2^{n,m}.
\end{split}\]

\paragraph{Step 1.} Assume that $n=1$. Fist, we show that:
\begin{equation} \label{expression of $L_{1,m}$}
L_{1,m}=F^{0,m}_1.
\end{equation}
The kernel of $F^{0,m}_2$ is equal to $\V{\otimes m-1}\otimes J_1=\V{\otimes m}$, that is, $F^{0,m}_2$ is equal to $0_{\V{\otimes m}}$. In particular, $F^{0,m}_1\vee F^{0,m}_2$ is equal to $F^{0,m}_1$. Moreover, the kernel of $F^{1,m}_1$ is equal to $I(R)_{m-1}\otimes V$ and the kernel of $F_1^{0,m}$ is equal to $I(R)_m$. The inclusion $I(R)_m\subset I(R)_{m-1}\otimes V$ implies that $F^{0,m}_1$ is smaller than $F^{1,m}_1$. Hence, Relation~\ref{expression of $L_{1,m}$} holds.

Assume that $m=N$. The kernel of $F_1^{1,N}$ is equal to $I(R)_{N-1}\otimes V=\{0\}$, that is, $F_1^{1,N}$ is equal to $\id{\V{\otimes N}}$. In particular, $R_{1,N}$ is equal to $F_2^{1,N}$. Moreover, we have:
\[\begin{split}
F_1^{0,N}&={\theta_{X^{(N)}}}^{-1}\left(I(R)_N\right)\\
&={\theta_{X^{(N)}}}^{-1}\left(R\right),\ \text{and}\\
F_2^{1,N}&={\theta_{X^{(N)}}}^{-1}\left(J_2\right)\\
&={\theta_{X^{(N)}}}^{-1}\left(R\right).
\end{split}\]
Thus $L_{1,N}$ and $R_{1,N}$ are equal.

Assume that $m\geq N+1$. From Lemma~\ref{appartence au treillis}, we have:
\[\begin{split}
F^{0,m}_1&=S^{(m)}_0\wedge\cdots\wedge S^{(m)}_{m-N},\\
F^{1,m}_1&=S^{(m)}_0\wedge\cdots\wedge S^{(m)}_{m-N-1},\\
F_2^{1,m}&=S^{(m)}_{m-N}.
\end{split}\]
Thus, $R_{1,m}$ is equal to $F^{0,m}_1$. We conclude that Proposition~\ref{utilisation $(*)$} holds for $n=1$ and $m\geq N$.

\paragraph{Step 2.} Assume that, $n\geq 2$ and that $l_N(n+1)\leq m < l_N(n+2)$. From Lemma~\ref{parity}, $m-l_N(n)$ is smaller than $N-1$. Thus, the kernel of $F^{n,m}_1$ is equal to $\{0\}$, that is, $F_1^{n,m}$ is equal to $\id{\V{\otimes m}}$. In particular, $L_{n,m}$ is equal to $F^{n-1,m}_1\vee F^{n-1,m}_2$ and $R_{n,m}$ is equal to $F^{n,m}_2$.

From Lemma~\ref{appartence au treillis}, we have:
\[\begin{split}
F^{n-1,m}_1&=S^{(m)}_0\wedge\cdots\wedge S^{(m)}_{m-l_N(n+1)},\\
F^{n-1,m}_2&=S^{(m)}_{m-l_N(n)}\vee\cdots\vee S^{(m)}_{m-N},\\
F^{n,m}_2&=S^{(m)}_{m-l_N(n+1)}\vee\cdots\vee S^{(m)}_{m-N}.
\end{split}\]
Moreover, from Lemma~\ref{lemme intermediaire 1}, we have:
\[\left(S^{(m)}_0\wedge\cdots\wedge S^{(m)}_{m-l_N(n+1)}\right)\vee S^{(m)}_{m-l_N(n)}=S^{(m)}_{m-l_N(n+1)}\vee\cdots\vee S^{(m)}_{m-l_N(n)}.\]
The law $\vee$ being associative, we deduce that Proposition~\ref{utilisation $(*)$} holds for every integers $n$ and $m$ such that $n\geq 2$ and that $l_N~(~n~+~1~)~\leq ~m~ < ~l_N~(~n~+~2~)$.

\paragraph{Step 3.} Assume that $n\geq 2$ and $m\geq l_N(n+2)$. From Lemma~\ref{appartence au treillis}, we have:
\[\begin{split}
F_1^{n-1,m}&=S^{(m)}_0\wedge\cdots\wedge S^{(m)}_{m-l_N(n+1)} ,\ \text{and}\\
F^{n,m}_1&=S^{(m)}_0\wedge\cdots\wedge S^{(m)}_{m-l_N(n+2)}.
\end{split}\]
Thus, letting $T_{n,m}=S^{(m)}_{m-l_N(n+2)+1}\wedge\cdots\wedge S^{(m)}_{m-l_N(n+1)}$, we have:
\[F^{n-1,m}_1=F^{n,m}_1\wedge T_{n,m}.\]
The lattice generated by $S^{(m)}_0,\cdots,S^{(m)}_{m-N}$ being distributive, we have:
\[F^{n-1,m}_1\vee F^{n-1,m}_2=\left(F^{n,m}_1\vee F^{n-1,m}_2\right)\wedge\left(T_{n,m}\vee F^{n-1,m}_2\right).\]
Using the inequality $F^{n,m}_1\preceq\left(F^{n,m}_1\vee F^{n-1,m}_2\right)$, we deduce:
\[L_{n,m}=F^{n,m}_1\wedge\left(T_{n,m}\vee F^{n-1,m}_2\right).\]
From Lemma~\ref{lemme intermediaire 2}, $T_{n,m}\vee F^{n-1,m}_2$ is equal to $F_2^{n,m}$. Thus, Proposition~\ref{utilisation $(*)$} holds for every integers $n$ and $m$ such that $n\geq 2$ and that $m\geq l_N(n+2)$.
\end{proof}

\subsubsection{Theorem}\label{Theorem}

\textit{Let $A$ be an $N$-homogeneous algebra. If $A$ admits an extra-confluent presentation $\EV{X\mid R}$, then the left bound of $\EV{X\mid R}$ is a contracting homotopy for the Koszul complex of $A$.}

\begin{proof}

Let $\phi$ be the endomorphism of $\tens{V}$ which maps any element to its unique normal form for $\EV{X\mid R}$.

The presentation $\EV{X\mid R}$ is side-confluent. Thus, from Proposition~\ref{relation r}, it is sufficient to show that for every integers $n$ and $m$ such that $n\geq 1$ and $m\geq l_N(n)$ we have:
\[\left(r_{n,m}\right)\ \ \ {F_1^{n,m}\wedge F_2^{n,m}}_{\mid K^{(m)}_n}={F_1^{n-1,m}\vee F_2^{n-1,m}}_{\mid K^{(m)}_n},\]
where $K^{(m)}_n$ is the vector space $\im{\phi_{\mid\V{\otimes m-l_N(n)}}}\otimes J_n$.\\

Assume that $l_N(n)\leq m<l_N(n+1)$. We show that $F_1^{n,m}\wedge F_2^{n,m}$ and $F_1^{n-1,m}\vee F_2^{n-1,m}$ are equal to $\id{\V{\otimes m}}$.

The hypothesis $l_N(n)\leq m<l_N(n+1)$ implies that $m-l_N(n)$ is smaller than $N-1$. In particular, the kernel of $F_1^{n,m}$ is equal to $\{0\}$, that is, $F_1^{n,m}$ is equal to $\id{\V{\otimes m}}$. Moreover, $F_2^{n,m}$ is also equal to $\id{\V{\otimes m}}$. Thus, the morphism $F_1^{n,m}\wedge F_2^{n,m}$ is equal to $\id{\V{\otimes m}}$. From Lemma~\ref{parity}, the morphism $F_1^{n-1,m}$ is equal to $\id{\V{\otimes m}}$. Thus $F_1^{n-1,m}\vee F_2^{n-1,m}$ is equal to $\id{\V{\otimes m}}$ and Relation $(r_{n,m})$ holds.\\

Assume that $m\geq l_N(n+1)$. From Lemma~\ref{commute} the operators $F^{n,m}_1$ and $F^{n-1,m}_1\vee F^{n-1,m}_2$ commute. We deduce from Relation~\ref{lower bound} (see page~\pageref{lower bound}):
\[F^{n,m}_1\wedge\left(F^{n-1,m}_1\vee F^{n-1,m}_2\right)=\left(F^{n-1,m}_1\vee F^{n-1,m}_2\right)F^{n,m}_1.\]
From Lemma~\ref{expression $F_1$}, the image of $F_1^{n,m}$ is equal to $\im{\phi_{\mid\V{\otimes m-l_N(n)}}}\otimes\V{\otimes l_N(n)}$. Thus, $K^{(m)}_n$ is included in $\im{F^{n,m}_1}$. Hence, the restriction of $F^{n,m}_1\wedge\left(F^{n-1,m}_1\vee F^{n-1,m}_2\right)$ to $K^{(m)}_n$ is equal to the restriction of $F_1^{n-1,m}\vee F_2^{n-1,m}$ to $K^{(m)}_n$. Moreover, the presentation $\EV{X\mid R}$ satisfies the extra-condition. Thus, from Proposition~\ref{utilisation $(*)$}, $F^{n,m}_1\wedge\left(F^{n-1,m}_1\vee F^{n-1,m}_2\right)$ is equal to $F_1^{n,m}\wedge F_2^{n,m}$. Hence, Relation $\left(r_{n,m}\right)$ holds.
\end{proof}

\section{Examples}\label{Examples}

In this section, we consider three examples of algebras which admit an extra-confluent presentation: the symmetric algebra, monomial algebras satisfying the overlap property and the Yang-Mills algebra over two generators. For each of these examples we explicit the left bound constructed in Section~\ref{Construction}. 

\subsection{The symmetric algebra}

In this section we consider the symmetric algebra $A=\K[x_1,\cdots ,x_d]$ over $d$ generators. This algebra admits the presentation $\EV{X\mid R}$ where $X$ is equal to $\left\{x_1,\cdots ,x_d\right\}$ and $R$ is equal to $\left\{x_ix_j=x_jx_i,\ 1\leq i\neq j\leq d\right\}$.

\subsubsection{Extra-confluence}

We consider the order $x_1<\cdots <x_d$. The operator $S\in\text{End}\left(\V{\otimes 2}\right)$ of the presentation $\EV{X\mid R}$ is defined on the basis $X^{(2)}$ by
\[S(x_ix_j)=\left\{\begin{split}
&x_jx_i,\ \text{if}\ i>j,\\
&x_ix_j,\ \text{otherwise}.
\end{split}\right.\]
Let $w=x_ix_jx_k\in X^{(3)}$. If $k$ is strictly smaller than $j$ and $i$ is strictly smaller than $k$, we have
\[\begin{split}
\EV{S\otimes\id{V},\id{V}\otimes S}^3(w)&=\EV{\id{V}\otimes S,S\otimes\id{V}}^3(w)\\
&=x_kx_jx_i.
\end{split}\]
In the other cases the elements $\EV{S\otimes\id{V},\id{V}\otimes S}^2(w)$ and $\EV{\id{V}\otimes S,S\otimes\id{V}}^2(w)$ are equal. In particular the two operators $\EV{S\otimes\id{V},\id{V}\otimes S}^3$ and $\EV{\id{V}\otimes S,S\otimes\id{V}}^3$ are equal.  Moreover, $N$ is equal to 2. Thus, from Remark~\ref{Extra-confluent presentation for quadratic algebras}, the presentation $\EV{X\mid R}$ is extra-confluent. The normal form of a word $x_{i_1}\cdots x_{i_n}$ is equal to $x_{j_1}\cdots x_{j_n}$ where $\{j_1,\cdots ,j_n\}=\{x_{i_1},\cdots ,x_{i_n}\}$ and $j_1\leq\cdots\leq j_n$.

\subsubsection{The Koszul complex of the symmetric algebra} 

The morphism $\partial_1:A\otimes V\F{}A$ is defined by $\partial_1(1_A\otimes x_i)=\overline{x_i}$, for every $1\leq i\leq d$. The morphism $\partial_2:A\otimes\overline{R}\F{}A\otimes V$ is defined by $\partial_2(1_A\otimes (x_jx_i-x_ix_j))=\overline{x_j}\otimes x_i-\overline{x_i}\otimes x_j$, for every $1\leq i<j\leq d$. If $d$ is greater than 3, the vector space $J_3$ is spanned by the elements
\[\begin{split}
r_{i_1<i_2<i_3}:&=x_{i_3}\left(x_{i_2}x_{i_1}-x_{i_1}x_{i_2}\right)-x_{i_2}\left(x_{i_3}x_{i_1}-x_{i_1}x_{i_3}\right)+x_{i_1}\left(x_{i_3}x_{i_2}-x_{i_2}x_{i_3}\right)\\
&=\left(x_{i_3}x_{i_2}-x_{i_2}x_{i_3}\right)x_{i_1}-\left(x_{i_3}x_{i_1}-x_{i_1}x_{i_3}\right)x_{i_2}+\left(x_{i_2}x_{i_1}-x_{i_1}x_{i_2}\right)x_{i_3},
\end{split}\]
where $1\leq i_1<i_2<i_3\leq d$. The morphism $\partial_3:A\otimes J_3\F{}A\otimes\overline{R}$ maps the element $1_A\otimes r_{i_1<i_2<i_3}$ to $\overline{x_{i_3}}\otimes\left(x_{i_2}x_{i_1}-x_{i_1}x_{i_2}\right)-\overline{x_{i_2}}\otimes\left(x_{i_3}x_{i_1}-x_{i_1}x_{i_3}\right)+\overline{x_{i_1}}\otimes\left(x_{i_3}x_{i_2}-x_{i_2}x_{i_3}\right)$.

Assume that $d$ is greater than 4 and let $n$ be an integer such that $3\leq n\leq d-1$. We denote by $I_n$ the set of sequences $i_1<\cdots <i_n$ such that $1\leq i_1$ and $i_n\leq d$. Assume that $r_l$ is defined for every $l\in I_{n}$. For every $l=i_1<\cdots <i_{n+1}\in I_{n+1}$ and every $1\leq j\leq n+1$ we denote by $l_j$ the element of $I_n$ obtained from $l$ removing $i_j$. Then, let
\[r_l=\sum_{j=0}^{n+1}(-1)^{-\eta(n+j)}x_{i_j}r_{l_j},\]
where $\eta:\N\F{}\{-1,1\}$ is defined by $\eta(k)=1$ if $k$ is even and $\eta(k)=-1$ if $k$ is odd. For every $4\leq n\leq d$, the vector space $J_n$ is spanned by the elements $r_l$ for $l\in I_n$. The map $\partial_{n}:A\otimes J_{n}\F{}A\otimes J_{n-1}$ is defined by
\[\partial_{n}\left(1_A\otimes r_{l}\right)=\sum_{j=1}^{n}(-1)^{-\eta(n-1+j)}\overline{x_{i_j}}\otimes r_{l_j}.\]
For every integer $n$ such that $n\geq d+1$, $J_n$ is equal to $\{0\}$. 

\subsubsection{The construction of $h_1$}

Let $m$ be an integer such that $m\geq 2$. Let $P_{1,m}~=~\left(F_1^{1,m},F_2^{1,m}\right)$ be the reduction pair of bi-degree $(1,m)$ associated with $\EV{X\mid R}$. The morphisms $F_1^{1,m}$ and $F_2^{1,m}$ are defined by
\[{F_1^{1,m}}\left(x_{i_1}\cdots x_{i_m}\right)=\widehat{w}x_{i_m},\ \text{where}\ w=x_{i_1}\cdots x_{i_{m-1}},\ \text{and}\]
\[F_2^{1,m}\left(x_{i_1}\cdots x_{i_m}\right)=x_{i_1}\cdots x_{i_{m-2}}\widehat{w},\ \text{where}\ w=x_{i_{m-1}}x_{i_m}.
\]
These morphisms satisfy the relation $\EV{F_1^{1,m},F_2^{1,m}}^4=\EV{F_2^{1,m},F_1^{1,m}}^3$. Thus, we consider the $P_{1,m}$-representation of $\M{A}_4$:
\[\begin{split}
  \varphi_{1,m}\colon\M{A}_{4}&\longrightarrow\text{End}\left(\V{\otimes m}\right),\\
  s_i&\longmapsto F_i^{1,m}
\end{split}
\]
The image of $\gamma_1=(1-s_2)(s_1+s_1s_2s_1)$ through this morphism is equal to $F_1^{1,m}-F_2^{1,m}F_1^{1,m}$. Let $wx_{i_1}\in X^{(m)}$. Denoting by $\widehat{w}=w'x_{i_2}$, $\varphi_{1,m}\left(\gamma_1\right)\left(wx_{i_1}\right)$ is equal to $w'\left(x_{i_2}x_{i_1}-x_{i_1}x_{i_2}\right)$ if $i_2<i_1$ and $\varphi_{1,m}\left(\gamma_1\right)\left(wx_{i_1}\right)$ is equal to $0$ otherwise. Then, the map $h_1:A\otimes V\F{} A\otimes\overline{R}$ is defined by
\[h_1\left(\overline{w}\otimes x_{i_1}\right)=\left\{\begin{split}
&\overline{w'}\otimes\left(x_{i_2}x_{i_1}-x_{i_1}x_{i_2}\right),\ \text{if}\  i_2<i_1,\\
&0,\ \text{otherwise}.
\end{split}\right.\]

\subsubsection{The construction of $h_2$}

Let $m$ be an integer such that $m\geq 3$. Let $P_{2,m}~=~\left(F_1^{2,m},F_2^{2,m}\right)$ be the reduction pair of bi-degree $(2,m)$ associated with $\EV{X\mid R}$. The morphisms $F_1^{2,m}$ and $F_2^{2,m}$ are defined by
\[{F_1^{1,m}}\left(x_{i_1}\cdots x_{i_m}\right)=\widehat{w}x_{i_{m-1}}x_{i_m},\ \text{where}\ w=x_{i_1}\cdots x_{i_{m-2}},\ \text{and}\]
\[F_2^{1,m}\left(x_{i_1}\cdots x_{i_m}\right)=\left\{\begin{split}
&x_{i_1}\cdots x_{i_m-3}\left(r_{i_{m-2}<i_{m-1}<i_m}\right),\ \text{if}\ i_{m-2}<i_{m-1}<i_m,\\
&0,\ \text{otherwise}.
\end{split}\right.
\]
These morphisms satisfy the relation $\EV{F_1^{2,m},F_2^{2,m}}^4=\EV{F_2^{2,m},F_1^{2,m}}^3$. Thus, we consider the $P_{2,m}$-representation of $\M{A}_4$:
\[\begin{split}
  \varphi_{2,m}\colon\M{A}_{4}&\longrightarrow\text{End}\left(\V{\otimes m}\right).\\
  s_i&\longmapsto F_i^{2,m}
\end{split}
\]
The image of $\gamma_1=(1-s_2)(s_1+s_1s_2s_1)$ through this morphism is equal to $F_1^{2,m}-F_2^{2,m}F_1^{2,m}$. Let $wx_{i_2}x_{i_1}\in X^{(m)}$. Denoting by $\widehat{w}=w'x_{i_3}$,  $\varphi_{2,m}\left(\gamma_1\right)\left(wx_{i_2}x_{i_1}\right)$ is equal to $w'  r_{i_1<i_2<i_3}$ if $i_1<i_2<i_3$ and $\varphi_{2,m}\left(\gamma_1\right)\left(wx_{i_2}x_{i_1}\right)$ is equal to $0$ otherwise. Then, the map $h_2:A\otimes\overline{R}\F{} A\otimes J_3$ is defined by
\[h_2\left(\overline{w}\otimes\left(x_{i_2}x_{i_1}-x_{i_1}x_{i_2}\right)\right)=\left\{\begin{split}
&\overline{w'}\otimes\left(r_{i_1<i_2<i_3}\right),\ \text{if}\  i_1<i_2<i_3,\\
&0,\ \text{otherwise}.
\end{split}\right.\]

\subsubsection{The construction of $h_n$}

More generally, for every $\overline{w}\otimes r_{i_1<\cdots <i_n}$ we denote by $\widehat{w}~=~w'x_{i_{n+1}}$. The map $h_n:A\otimes J_n\F{} A\otimes J_{n+1}$ is defined by
\[h_n\left(\overline{w}\otimes r_{i_1<\cdots <i_n}\right)=\left\{\begin{split}
&\overline{w'}\otimes r_{i_1<\cdots <i_{n+1}},\ \text{if}\  i_{1}<\cdots <i_{n+1},\\
&0,\ \text{otherwise}.
\end{split}\right.\]

\subsubsection{Remark}

The left bound family of $\EV{X\mid R}$ is the contracting homotopy constructed in the proof of~\cite[Proposition 3.4.13]{MR2954392}.

\subsection{Monomial algebras satisfying the overlap property}

In the section we consider the example from~\cite[Proposition 3.8]{MR1832913}. We consider a \emph{monomial algebra} $A$ over $d$ generators: $X=\{x_1,\cdots ,x_d\}$ and $R=\{w_1,\cdots ,w_l\}$ is a set of words of length $N$. We assume that the presentation $\EV{X\mid R}$ satisfies the \emph{overlap property}. This property is stated as follows: 

\subsubsection{The overlap property}

For every integer $n$ such that $N+2\leq n\leq 2N-1$ and for any word $w=x_{i_1}\cdots x_{i_n}$ such that $x_{i_1}\cdots x_{i_N}$ and $x_{i_{n-N+1}}\cdots x_{i_n}$ belong to $R$, all the sub-words of length $N$ of $w$ belong to $R$.

\subsubsection{Extra-confluence}

For any choice of order on $X$, the operator $S\in\text{End}\left(\V{\otimes N}\right)$ of the presentation $\EV{X\mid R}$ is defined on the basis $X^{(N)}$ by
\[S(w)=\left\{\begin{split}
&0,\ \text{if}\ w\in R,\\
&w,\ \text{otherwise}.
\end{split}\right.\]
As a consequence, for every integer $m$ such that $1\leq m\leq N-1$, the operators $S\otimes\id{\V{\otimes m}}$ and $\id{\V{\otimes m}}\otimes S$ commute. Thus, the presentation $\EV{X\mid R}$ is side-confluent. Moreover, for monomial algebras, the extra-condition is equivalent to the overlap property. Thus, the presentation $\EV{X\mid R}$ is extra-confluent. The normal form of a word $w$ is equal to $0$ if $w$ admits a sub-word which belongs to $R$, and $w$ otherwise. 

\subsubsection{The Koszul complex of a monomial algebra}

Let $n$ be an integer such that $n\geq 2$. The vector space $J_n$ is spanned by words $w$ of length $l_N(n)$ such that every sub-word of length $N$ of $w$ belongs to $R$. The morphism $\partial_n:A\otimes J_n\F{} A\otimes J_{n-1}$ maps $1_A\otimes x_{i_1}\cdots x_{i_{l_N(n)}}$ to $\overline{w'}\otimes x_{i_{l_N(n)-l_N(n-1)+1}}\cdots x_{i_{l_N(n)}}$, where $w'$ is equal to $x_{i_1}\cdots x_{i_{l_N(n)-l_N(n-1)}}$.

\subsubsection{The contracting homotopy}

Let $n$ and $m$ be two integers such that $m\geq l_n(n)$. Let $P_{n,m}~=~\left(F_1^{n,m},F_2^{n,m}\right)$ be the reduction pair of bi-degree $(n,m)$ associated with $\EV{X\mid R}$. The operators $F_1^{n,m}$ and $F_2^{n,m}$ are defined by
\[F_1^{n,m}(x_{i_1}\cdots x_{i_m})=\left\{\begin{split}
&0,\ \text{if a sub-word of length}\ N\ \text{of}\ x_{i_1}\cdots x_{i_{m-l_n(n)}}\ \text{belongs to}\ R,\\
&w,\ \text{otherwise},
\end{split}\right.\]
and 
\[F_2^{n,m}(x_{i_1}\cdots x_{i_m})=\left\{\begin{split}
&0,\ \text{if}\ x_{i_{m-l_N(n+1)+1}}\cdots x_{i_m}\in J_{n+1},\\
&w,\ \text{otherwise}.
\end{split}\right.\]
These operators commute. Thus, we consider the $P_{n,m}$-representation of $\M{A}_2$:
\[\begin{split}
  \varphi_{n,m}\colon\M{A}_{2}&\longrightarrow\text{End}\left(\V{\otimes m}\right).\\
  s_i&\longmapsto F_i^{2,m}
\end{split}
\]
The image of $\gamma_1=(1-s_2)s_1$ through this morphism is equal to $F_1^{n,m}-F_2^{n,m}F_1^{n,m}$. Let $w~=~x_{i_1}\cdots x_{i_m}$ be an element of $X^{(m)}$. If $w$ is such that no sub-word of length $N$ of $x_{i_1}\cdots x_{i_{m-l_n(n)}}$ belongs to $R$ and if $x_{i_{m-l_N(n+1)+1}}\cdots x_{i_m}$ belongs to $J_{n+1}$, $\varphi_{n,m}\left(\gamma_1\right)(w)$ is equal to $w$. In the other cases $\varphi_{n,m}\left(\gamma_1\right)(w)$ is equal to $0$. Then, the morphism $h_n:A\otimes J_n\F{} A\otimes J_{n+1}$ is defined by
\[h_n\left(\overline{w}\otimes x_{i_{m-l_N(n)+1}}\cdots x_{i_m}\right)=\left\{\begin{split}
& \overline{w'}\otimes x_{i_{m-l_N(n+1)+1}}\cdots x_{i_m},\ \text{if}\ x_{i_{m-l_N(n+1)+1}}\cdots x_{i_m}\in J_{n+1}, \\
&0,\ \text{otherwise},
\end{split}
\right.\]
where $w=x_{i_1}\cdots x_{i_{m-l_N(n)}}$ and $w'=x_{i_1}\cdots x_{i_{m-l_N(n+1)}}$.

\subsection{The Yang-Mills algebra over two generators}\label{The Yang-Mills algebra over two generators}

In this section we explicit the left upper bound family associated with the presentation $\EV{X\mid R}$ from Example~\ref{Yang-Mills side-confluent} of the Yang-Mills algebra over two generators. Recall that $X=\{x_1,x_2\}$ and $R=\{f_1,f_2\}$ where
\[\begin{split}
f_1&=x_2x_1x_1-2x_1x_2x_1+x_1x_1x_2,\ \text{and} \\
f_2&=x_2x_2x_1-2x_2x_1x_2+x_1x_2x_2 .
\end{split}\] 
The acyclicty of the Koszul complex of this algebra was proven in~\cite[Section 6.3]{MR3299599} using the arguments of Example~\ref{Yang-Mills extra-confluent}. In this section, we propose an other proof, based on the construction of an explicit contracting homotopy.

\subsubsection{Extra-confluence}

Recall that for the order $x_1<x_2$, the operator $S\in\text{End}\left(\V{\otimes 3}\right)$ of the presentation $\EV{X\mid R}$ is defined on the basis $X^{(3)}$ by
\[S(w)=\left\{\begin{split}
&2x_1x_2x_1-x_1x_1x_2,\ \text{if}\ w=x_2x_1x_1, \\
&2x_2x_1x_2-x_1x_2x_2,\ \text{if}\ w=x_2x_2x_1, \\
&w,\ \text{otherwise}.
\end{split}\right.\]
Recall from Example~\ref{Yang-Mills extra-confluent} that this presentation is extra-confluent.

\subsubsection{The Koszul complex of the Yang-Mills algebra}

The morphism $\partial_1:A\otimes V\F{} A$ is defined by $\partial_1\left(1_A\otimes x_i\right)=\overline{x_i}$ for $i=1$ or $2$. The morphism $\partial_2:A\otimes\overline{R}\F{} A\otimes V$ is defined by \[\begin{split}
\partial_2\left(1_A\otimes f_1\right)&=\overline{x_2x_1}\otimes x_1-2\overline{x_1x_2}\otimes x_1+\overline{x_1x_1}\otimes x_2,\ \text{and}\\
\partial_2\left(1_A\otimes f_2\right)&=\overline{x_2x_2}\otimes x_1-2\overline{x_2x_1}\otimes x_2+\overline{x_1x_2}\otimes x_2.
\end{split}\]
The vector space $J_3=\left(V\otimes\overline{R}\right)\cap\left(\overline{R}\otimes V\right)$ is the one-dimensional vector space spanned by
\[\begin{split}
v&=x_2f_1+x_1f_2\\
&=f_2x_1+f_1x_2.
\end{split}\]
The morphism $\partial_3:A\otimes J_3\F{} A\otimes\overline{R}$ is defined by
\[\partial_3\left(1_A\otimes v\right)=\overline{x_2}\otimes f_1+\overline{x_1}\otimes f_2.\]
For every integer $n$ such that $n\geq 4$, the vector space $J_n$ is equal to $\{0\}$.

\subsubsection{The construction of $h_1$}

Recall from Proposition~\ref{Diamond lemma} that the algebra $A$ admits as a basis the set $\left\{\overline{w},\ w\in\EV{X}\ \text{is a normal form}\right\}$. Thus, it is sufficient to define $h_1\left(\overline{w}\otimes x_i\right)$ for every normal form $w\in\EV{X}$ and $i=1$ or $2$.\\

Let $m$ be an integer such that $m\geq 3$. Let $P_{1,m}=\left(F_1^{1,m},F_2^{1,m}\right)$ be the reduction pair of bi-degree $(1,m)$ associated with $\EV{X\mid R}$. The morphisms $F_1^{1,m}$ and $F_2^{1,m}$ are defined by
\[{F_1^{1,m}}\left(x_{i_1}\cdots x_{i_m}\right)=\widehat{w}x_{i_m},\ \text{where}\ w=x_{i_1}\cdots x_{i_{m-1}},\ \text{and}\]
\[F_2^{1,m}\left(x_{i_1}\cdots x_{i_m}\right)=x_{i_1}\cdots x_{i_{m-3}}\widehat{w},\ \text{where}\ w=x_{i_{m-2}}x_{i_{m-1}}x_{i_m}.
\]
These morphisms commute. Thus, we consider the $P_{1,m}$-representation of $\M{A}_2$:
\[\begin{split}
  \varphi_{1,m}\colon\M{A}_{2}&\longrightarrow\text{End}\left(\V{\otimes m}\right).\\
  s_i&\longmapsto F_i^{1,m}
\end{split}
\]
The image of $\gamma_1=(1-s_2)s_1$ through this morphism is equal to $F_1^{1,m}-F_2^{1,m}F_1^{1,m}$.

Let $w$ be a normal form such that the length of $w$ is equal to $m- 1$. The word $wx_2$ does not factorize on the right by $x_2x_1x_1$ or $x_2x_2x_1$. Thus, $\varphi_{1,m}\left(\gamma_1\right)(wx_2)$ is equal to $0$. In particular, $h_1\left(\overline{w}\otimes x_2\right)$ is equal to $0$ for every normal form $w\in\EV{X}$. If $w$ does not factorize on the right by $x_2x_1$ or $x_2x_2$, $\varphi_{1,m}\left(\gamma_1\right)(wx_1)$ is equal to $0$. Thus, $h_1\left(\overline{w}\otimes x_1\right)$ is equal to $0$ for every normal form $w\in\EV{X}$ such that $w$ does not factorize on the right by $x_2x_1$ or $x_2x_2$. If $w$ can be written $w'x_2x_1$ (respectively $w'x_2x_2$), then $\varphi_{1,m}\left(\gamma_1\right)(wx_1)$ is equal to $w'\left(2x_1x_2x_1-x_1x_1x_2\right)$ (respectively $w'\left(2x_2x_1x_2-x_1x_2x_2\right)$). Thus, we have:
\[h_1\left(\overline{w}\otimes x_1\right)=\left\{
\begin{split}
& \overline{w'}\otimes\left(2x_1x_2x_1-x_1x_1x_2\right),\ \text{if}\ w=w'x_2x_1,\\
& \overline{w'}\otimes\left(2x_2x_1x_2-x_1x_2x_2\right),\ \text{if}\ w=w'x_2x_2.
\end{split}
\right.\]

\subsubsection{The construction of $h_2$}

Recall from Proposition~\ref{Diamond lemma} that the algebra $A$ admits as a basis the set $\left\{\overline{w},\ w\in\EV{X}\ \text{is a normal form}\right\}$. Thus, it is sufficient to define $h_2\left(\overline{w}\otimes f_i\right)$ for every normal form $w\in\EV{X}$ and $i=1$ or $2$.\\

Let $m$ be an integer such that $m\geq 4$. Let $P_{2,m}=\left(F_1^{2,m},F_2^{2,m}\right)$ be the reduction pair of bi-degree $(2,m)$ associated with $\EV{X\mid R}$. The operator $F_1^{2,m}$ maps a word $w\in X^{(m)}$ to $\widehat{w_1}w_2$, where $w_1\in\EV{X}$ and $w_2\in X^{(4)}$ are such that $w=w_1w_2$. The operator $F_2^{2,m}$ is equal to $\id{\V{\otimes m-4}}\otimes F$ where $F$ is equal to $\theta_{X^{(4)}}^{-1}\left(J_3\right)$. The kernel of $F$ is the one-dimensional vector space spanned by $v$. Thus, $F(\lm{v})$ is equal to $\lm{v}-v$, and for every $w\in X^{(4)}\setminus\{\lm{v}\}$, $F(w)$ is equal to $w$. Thus, $F$ is defined on the basis $X^{(4)}$ by
\[F(w)=\left\{
\begin{split}
&2x_2x_1x_2x_1-x_2x_1x_1x_2-x_1x_2x_2x_1+2x_1x_2x_1x_2-x_1x_1x_2x_2,\ \text{if}\ w=x_2x_2x_1x_1,\\
&w,\ \text{otherwise}.
\end{split}
\right.\]

The two operators $F_1^{2,m}$ and $F_2^{2,m}$ commute. Thus, we consider the $P_{2,m}$-representation of $\M{A}_2$:
\[\begin{split}
  \varphi_{2,m}\colon\M{A}_{2}&\longrightarrow\text{End}\left(\V{\otimes m}\right).\\
  s_i&\longmapsto F_i^{2,m}
\end{split}
\]
The image of $\gamma_1=(1-s_2)s_1$ is equal to $F_1^{2,m}-F_2^{2,m}F_1^{2,m}$.

Let $w$ be a normal form such that the length of $w$ is equal to $m- 1$. The word $x_2x_2x_1x_1$ does not occur in the decomposition of $wf_2$. Thus, $\varphi_{2,m}(wf_2)$ is equal to $0$. In particular $h_2\left(\overline{w}\otimes f_2\right)$ is equal to $0$ for every normal form $w\in\EV{X}$. If $w$ does not factorize on the right by $x_2$, the word $x_2x_2x_1x_1$ does not occur in the decomposition of $wf_1$. Thus, $\varphi_{2,m}(wf_1)$ is equal to $0$. In particular $h_2\left(\overline{w}\otimes f_1\right)$ is equal to $0$ for every normal form $w\in\EV{X}$ such that $w$ does not factorize on the right by $x_2$. Assume that $w$ factorize on the right by $x_2$: $w=w'x_2$. Thus, $\varphi_{2,m}(wf_1)$ is equal to $w'\left(x_2x_2x_1x_1-F(x_2x_2x_1x_1)\right)$. In this case we have
\[h_2\left(\overline{w}\otimes f_1\right)=\overline{w'}\otimes\left( x_2f_1+x_1f_2\right).\]

\bibliography{Biblio}

\newpage\noindent
\textbf{Cyrille Chenavier}\\
\textbf{INRIA, \'equipe $\pi r^2$}\\
\textbf{Laboratoire Preuves, Programmes et Syst\`emes, CNRS UMR 7126}\\
\textbf{Universit\'e Paris-Diderot}\\
\textbf{Case 7014}\\
\textbf{75205 PARIS Cedex 13}\\
\textbf{cyrille.chenavier@pps.univ-paris-diderot.fr}

\end{document}